 \documentclass[10pt]{article}
 \usepackage[top=1in,bottom=1in,left=1.5in,right=1.5in]{geometry}
  \usepackage{amsmath,amssymb}
  \usepackage{latexsym}
 \usepackage[dvips]{pstricks} 
  \usepackage{pst-node}
  \usepackage[dvips]{graphicx}

  \newcommand{\C}{\mathbb{C}}
  
  \newcommand{\N}{\mathbb{N}}
  \renewcommand{\P}{\mathbb{P}}

  \renewcommand{\a}{\mathbf{a}}

  \renewcommand{\c}{\mathbf{c}}
  \newcommand{\e}{\mathbf{e}}
  \newcommand{\f}{\mathbf{f}}
  \newcommand{\g}{\mathbf{g}}
  \newcommand{\gl}{\mathbf{GL}}

  \newcommand{\bS}{\mathbf{S}}
  
  \newcommand{\D}{\mathbf{D}}
  \newcommand{\T}{\mathbf{T}}
  \newcommand{\U}{\mathbf{U}}
  \renewcommand{\u}{\mathbf{u}}
  \renewcommand{\v}{\mathbf{v}}
  \newcommand{\V}{\mathbf{V}}
  \newcommand{\w}{\mathbf{w}}
  \newcommand{\W}{\mathbf{W}}
  \newcommand{\x}{\mathbf{x}}
  \newcommand{\X}{\mathbf{X}}
  \newcommand{\y}{\mathbf{y}}
  \newcommand{\Y}{\mathbf{Y}}
  \newcommand{\z}{\mathbf{z}}
  \newcommand{\0}{\mathbf{0}}
  \newcommand{\1}{\mathbf{1}}

  \newcommand{\cS}{\mathcal{S}}
  \newcommand{\cT}{\mathcal{T}}

  \newcommand{\lan}{\langle}
  \newcommand{\ran}{\rangle}
  \newcommand{\an}[1]{\lan#1\ran}
  \def\diag{\mathop{{\rm diag}}\nolimits}
  \newcommand{\hs}{\hspace*{\parindent}}
  \newcommand{\proof}{\hs \textbf{Proof.\ }}
  
  \newcommand{\tr}{\mathop{\mathrm{tr}}\nolimits}

  \newcommand{\Gr}{\mathop{\mathrm{Gr}}\nolimits}
  \newcommand{\trans}{^\top}
  \newcommand{\qed}{\hspace*{\fill} $\Box$\\}

  \newcommand{\codim}{\mathrm{codim}\;}

  \newcommand{\grank}{\mathrm{grank}}
  
  \newcommand{\maxrank}{\mathrm{mrank}\;}

  \newcommand{\rC}{\mathrm{C}}

  \newcommand{\rS}{\mathrm{S}}

  \newcommand{\adj}{\mathrm{adj\;}}
  \newcommand{\rank}{\mathrm{rank\;}}
  \newcommand{\brank}{\mathrm{brank\;}}

  \newcommand{\trace}{\mathrm{trace\;}}

  \newtheorem{theo}{\bfseries \hs Theorem}[section]
  
  \newtheorem{prop}[theo]{\bfseries \hs Proposition}
  
  \newtheorem{lemma}[theo]{\bfseries \hs Lemma}
  \newtheorem{corol}[theo]{\bfseries \hs Corollary}

  \numberwithin{equation}{section} 

 \renewcommand{\span}{\mathrm{span}}

\begin{document}

 \title{On tensors of border rank $l$ in $\C^{m\times n\times l}$}
 \author{
 Shmuel Friedland\\
 Department of Mathematics, Statistics and Computer Science\\
 University of Illinois at Chicago\\ Chicago, Illinois 60607-7045,
 USA\\ \texttt{E-mail: friedlan@uic.edu}
 }

 \date{November 12, 2010}
 \maketitle
 \begin{abstract}
 We study tensors in $\C^{m\times n\times l}$ whose border rank is $l$.
 We give a set-theoretic characterization of tensors in $\C^{3\times 3\times 4}$ and in $\C^{4\times 4\times 4}$ of border rank $4$ at most.
 \end{abstract}



 \section{Introduction}  Denote by $\C^{m\times n},\rS(m,\C),\C^{m\times n\times l}$ the linear spaces of $m\times n$ matrices, $m\times m$ symmetric matrices and $3$-tensors  $\cT=[t_{i,j,k}]_{i=j=k=1}^{m,n,l}$ of dimension $m\times
 n\times l$ over the field of complex numbers $\C$ respectively.  We identify $\C^{m\times n\times l}$ with
 $\C^{m}\otimes \C^{n}\otimes \C^l$.
 A rank one tensor is $\u\otimes\v\otimes\w=[u_iv_jw_k]\in \C^{m\times n\times l}$,
 where $\u=(u_1,\ldots,u_m)\trans\in\C^m,\v=(v_1,\ldots,v_n)\trans\in\C^n,\w=(w_1,\ldots,w_l)\trans\in\C^l$,
 and $\u,\v,\w$ are nonzero vectors.  The rank of a nonzero tensor $\cT$ is the minimal number $r:=\rank\cT$, such that $\cT=\sum_{i=1}^r \u_i\otimes\v_i\otimes\w_i$.  The \emph{border} rank
 of $\cT\ne 0$, denoted by $\brank \cT$, is a positive integer $q$, such that the following conditions
 hold.  First, $\cT$ is a limit of a sequence of tensors $\cT_{\nu}\in\C^{m\times n\times l},
 \rank \cT_{\nu}=q,\nu\in\N$.  Second, $\cT$ is not a limit of any sequence of tensors, such that each tensor in the sequence has rank $q-1$ at most.
 For $m,n,l\ge 2$ there exists tensors with $\rank \cT> \brank \cT$.  (This inequality does not hold for
 matrices, i.e. $\C^{m\times n}$.)  The maximal border rank of $\cT\in \C^{m\times n\times l}$ is called
 the generic rank of $\C^{m\times n\times l}$ and is denoted by $\grank(m,n,l)$.
 The value of $\grank(m,n,l)$ is known for many triples of $(m,n,l)$.
 The conjectural value of $\grank(m,n,l)$ is given in \cite{Fri08}.
 For $r \le \grank(m,n,l)$ denote by $V_r(m,n,l)\subset \C^{m\times n\times l}$
 the set of all $3$-tensors of border rank $r$ at most.  It is easy to see that $V_r(m,n,l)$
 is an irreducible variety in $\C^{m\times n\times l}$, which is a zero set of a number of homogeneous polynomials.  In fact, its projectivization
 is the $r$-secant variety of $\P^{m-1}\times \P^{n-1}\times \P^{l-1}$.

 A general problem is to characterize $V_{r}(m,n,l)$ in one of the following ways: set-theoretic,
 scheme theoretic and ideal theoretic.
 An elegant result of Strassen characterizes $V_{4}(3,3,3)$ \cite{Str}.  It is a hypersurface given
 by a homogeneous polynomial of degree $9$.
 This paper solves the set-theoretic aspect of the \textbf{Challenge Problem} posed by Elizabeth S. Allman
 in March 2007 (http://www.dms.uaf.edu/$\sim$eallman/): Determine the ideal defining the fourth secant variety of $\P^3 \times \P^3 \times \P^3$.  See \cite{AR07} for more details how this particular problem related to phylogenic ideals and varieties for general Markov models.
 The \emph{Salmon conjecture}
 \cite[Conjecture 3.24]{PS} stated that $V_{4}(4,4,4)$ is defined by polynomials of degree $5$ and $9$.
 A first nontrivial step in characterizing $V_{4}(4,4,4)$ is to characterize $V_4(3,3,4)$.
 It is shown in \cite{LM04} that $V_4(3,3,4)$ satisfies certain polynomial equations of degree $6$.
 (See also \cite[Remark 5.7]{LM06} and \cite{OB10}.)
 Hence the corrected version of the Salmon conjecture
 states that $V_{4}(4,4,4)$ is defined by polynomials of degree $5$, $6$ and $9$
 \cite[\S2]{Stu08}.

 The first main result of this paper shows that $V_{4}(3,3,4)$ is cut out by a set of polynomials
 of degree $9$ and $16$.
 Our second main result shows that $V_{4}(4,4,4)$ is cut out by a set of polynomials
 of degree $5$, $9$ and $16$.
 Most of the results in this paper are derived
 from results in matrix theory and relatively basic results in algebraic geometry.
 Whenever we could, we  stated our results in a general setting.

 We first explain briefly the main steps in the set-theoretic characterization of $V_4(4,4,4)$.
 First observe that $\cT\in\C^{4\times 4\times 4}$ is given by any four $p$-slices of $\cT$, for $p=1,2,3$. For example, the $(i,1)$-slice of $\cT$ is $T_{i,1}=[t_{i,j,k}]_{j,k=1}^4\in\C^{4\times 4}$ for $i=1,2,3,4$.  Let $\T_1(\cT)=\span(T_{1,1},\ldots,T_{4,1})$.
 Assume first the generic case that $\dim\T_1(\cT)=4$ and $\rank \cT=4$.   \cite[Theorem 2.4]{Fri08}
 yields that $\T_1(\cT)=\span(\u_1\v_1\trans,\ldots,\u_4\v_4\trans)$.  Assume the generic case that $\u_1,\ldots,\u_4\in\C^4$ and $\v_1,\ldots,\v_4\in\C^4$ are linearly independent.
 Let $Q\in\gl(4,\C)$ satisfy $Q\u_i=\w_i, \v_i\trans \w_j=\delta_{ij}, i,j=1,\ldots,4$.  Then any two matrices in $Q\T_1(\cT)$ commute.  (This result is well known to the experts, e.g. \cite{DeL06} and references
 therein.)
 This result is equivalent to the statement that for any $X,Y,Z\in\T_1(\cT)$ the following condition holds.
 \begin{equation}\label{XYZcomutcon}
 X(\adj Y) Z-Z(\adj Y) X=0.
 \end{equation}
 ($\adj Y$ is the adjoint matrix of $Y$.)
 These identities give rise to a system of homogeneous equations of degrees $5$ in the entries of $X,Y,Z$, which always hold if $\brank\cT\le 4$.
 Vice versa, if the above equalities hold and $\T_1(\cT)$ contains an invertible matrix then the results
 in \cite{GS} yields that $\brank\cT\le 4$.

 We next consider the case where $\T_1(\cT),\T_2(\cT),\T_3(\cT)$ does not contain an invertible matrix, and
 each three matrices in $\T_i(\cT)$ satisfies (\ref{XYZcomutcon}) for $i=1,2,3$.
 In \S5 we show that either $\brank\cT\le 4$ or by permuting factors in $\C^4\otimes\C^4\otimes \C^4$, if necessary, and changing bases in the first two components of $\C^4\times\C^4\times\C^4$
 $\cT$ can be viewed as a tensor $\C^{3\times 3\times 4}$.  This is Corollary 5.6 of \cite{LM06}.
 In \cite{LM06} this corollary is deduced from \cite[Prop. 5.4]{LM06}. Unfortunately, Proposition 5.4 is wrong, see \S5.

 Assume that $\cT\in V_4(3,3,4)$ and $\rank\cT=4$.  Then $\T_3(\cT)\subset \C^{3\times 3}$ is spanned by $\u_i\v_i\trans$ for $i=1,\ldots,4$.  Assume the generic case where any three vectors out of $\u_1,\ldots,\u_4$ and out of $\v_1,\ldots,\v_4$ are linearly independent.  Then there exists $L,R\in\gl(3,\C)$ such that $L\T_3(\cT)$ and $\T_3(\cT) R$ are $4$ dimensional subspaces of the subspace of $3\times 3$ symmetric matrices $\rS(3,\C)$, which are spanned by $4$ rank one symmetric matrices.
 Furthermore $L,R$ are unique up to a multiplication by a nonzero scalar, and
 \begin{equation}\label{RLident}
 LR\trans=R\trans L=\frac{\tr(LR\trans)}{3}I_3.
 \end{equation}
 The existence of nonzero $L,R$ such that $L\T_3(\cT),\T_3(\cT) R\subset \rS(3,\C)$
 is equivalent to the condition that the corresponding system of homogeneous linear equations in the entries of $L$ and $R$ respectively, given respectively by the coefficient matrices $C_L(\cT),C_R(\cT)\in \C^{12\times 9}$, have nontrivial solution. (Note that $L$ and $R$ have $9$ entries.)
 The entries of $C_L(\cT)$ and $C_R(\cT)$ are linear combinations of the entries of $\cT$ with coefficients $0,1,-1$. A necessary and sufficient condition for a nontrivial solution $R$ and $L$ is that all $9\times 9$ minors of $C_L(\cT)$ and $C_R(\cT)$ are zero.  These gives rise to a number of polynomial equations of degree $9$ that the entries of $\cT\in\C^{3\times 3\times 4}$ must satisfy.
 We show that the Strassen condition corresponds to some of the above polynomial equations.
 If the ranks of $C_L(\cT)$ and $C_R(\cT)$ is $8$ then $L$ and $R$ are determined uniquely up to a multiplication of a nonzero scalar.  The entries of $L$ and $R$ are polynomials of degree $8$ in the entries of $C_L(\cT)$ and $C_R(\cT)$.  The condition (\ref{RLident}) translates to a system of polynomial
 equations of degree $16$ in the entries of $\cT$.

 Assume first that $L$ is invertible.  Then $L\T_3(\cT)$ is a $4$-dimensional subspace of $\rS(3,\C)$.
 A $4$-dimensional generic subspace of $\rS(3,\C)$ spanned by $4$ rank one matrices.  Hence the assumption
 that $L\in\gl(3,\C)$ yields that $\cT\in V_4(3,3,4)$.  (This case does not require (\ref{RLident}).)  Assume that neither $L$ or $R$ are invertible.  Then
 $LR=0$.  In that case we also show that the above conditions imply that $\cT\in V_4(3,3,4)$, by considering a few cases.  (We need (\ref{RLident}) to rule out certain cases.)

 We survey briefly the contents of the paper.  In \S2 we discuss some known results which are needed in the
 next sections.  We recall a simple known condition that $\rank \cT$ is the dimension of the minimal
 subspace spanned by rank one matrices that contains the subspace spanned by $p$-slices of $\cT$, denoted by $\T_p(\cT)$, for each $p=1,2,3$.  Next we discuss a simple dimension condition on a generic subspace in $\U\subset\C^{m\times n}$ which implies that $\U$ is spanned by rank one matrices.  We translate this result to the border rank of $\cT$.  We recall the Strassen characterization of $V_4(3,3,3)$.
 Next we show that for a generic tensor $\cT\in\C^{m\times m\times l}$ of rank $m$ one can change a basis
 in the first factor of $\C^m$ such that the $l$ $3$-slices of $\cT$ are commuting matrices.
 These conditions give rise to the equations of type (\ref{XYZcomutcon}).  In \S3 we
 characterize subspaces $\U\subset\C^{m\times m}$ such that any $3$ matrices satisfy the condition (\ref{XYZcomutcon}) and most of the matrices in $\U$ have rank $m-1$.
 In \S4 we characterize $V_4(3,3,4)$.  \S5 we discuss necessary and sufficient
 conditions for a $\cT\in\C^{4\times 4\times 4}$ to have border rank $4$ at most.
 The nontrivial case is when $\T_1(\cT),\T_2(\cT),\T_3(\cT)$ do not contain an invertible matrix,
 and $\cT$ can not be viewed as a tensor in $\C^{3\times 3\times 4}$.
 We analyze the case where most of the matrices in $\T_3(\cT)$ have rank $2$.
 In this case the condition (\ref{XYZcomutcon}) holds trivially for $\T_3(\cT)$.
 We show that in some cases $\brank\cT=5$ if we do not assume the conditions (\ref{XYZcomutcon})
 for $\T_1(\cT)$.  This gives a counter-example to \cite[Proposition 5.4]{LM06}, and invalidates
 the proof of \cite[Corollary 5.6]{LM06}.  On the other hand we show that if (\ref{XYZcomutcon})
 holds for $\T_1(\cT)$ and of the matrices in $\T_3(\cT)$ have rank $2$, then $\brank \cT\le 4$.
 Most of $\S5$ is devoted to the proof of \cite[Corollary 5.6]{LM06}.  We need to consider the case
 where most of the matrices in $\T_p(\cT)$ have rank $3$ for $p=1,2,3$.
 Our analysis depends on the results in \S3.

 In \S6 we outline how to translate the
 problem of classifying tensors in $\C^{m\times n\times l}$ of rank $l$ if either $2\le l\le m,n$ or $m=n=l-1$ and $l\ge 4$.  It turns out that in the generic case this problem boils down to the condition that
 a corresponding $l$ subspace denote by $\bS(\T_3(\cT))\subset \rS(l,\C)$ is congruent to a subspace
 spanned by $l$ diagonal matrices.
 Note that the a simultaneous matrix diagonalization by congruence arises naturally in finding the rank
 decomposition of tensors \cite{DeL06}.
 We point out how some of these results can be generalized to tensors
 of border rank $l$ at most.

 \section{Preliminary results}

 We first recall a basic result on the rank of $3$ tensor $\cT=[t_{i,j,k}]_{i=j=k}^{m,n,l}\in\C^{n\times m\times l} $
 which is well known to the experts.  (See for example \cite[Theorem 2.4]{Fri08}.)
 By a $(k,3)$-slice, we denote the matrix
 $T_{k,3}(\cT)=T_{k,3}:=[t_{i,j,k}]_{i=j=1}^{m,n}\in\C^{m\times n}$ for $k=1,\ldots,l$.  Let $\T_3(\cT)
 :=\span(T_{1,3},\ldots,T_{l,3}) \subset \C^{m\times n}$.  We call $\T_3(\cT)$ the $3-rd$ subspace induced
 by $\cT$.
 \begin{theo}\label{bascharrnk} Let $\cT\in\C^{m\times n \times l}$.  Then $\rank \cT$ is the minimal
 dimension of a subspace $\U\subset \C^{m\times n}$ that contains $\T_3(\cT)$ and is spanned by rank
 one matrices.
 \end{theo}

 We can define similarly $(p,1)$ and $(q,2)$ slices of $\cT$ and the corresponding subspaces
 $\T_1(\cT),\T_2(\cT)$.  Hence Theorem \ref{bascharrnk} can be stated for $\T_1(\cT)$ and $\T_2(\cT)$
 respectively.  Also note that the space $\C^{m\times n\times l}$ can be identified with
 $\C^{m}\otimes \C^{n}\otimes \C^{l}$, and hence if we permute the three tensor factors $\C^m,\C^n,\C^l$
 we obtain $6$ isomorphic spaces of tensors.

 Let
 $$P=[p_{i'i}]\in\gl(m,\C),\; Q=[q_{j'j}]\in\gl(n,\C),\; R=[r_{k'k}]\in\gl(l,\C).$$
 Then we can change
 the standard bases
 in $\C^m,\C^n,\C^l$ to the bases obtained from the columns of $P^{-1},Q^{-1},R^{-1}$ respectively.
 In the new bases $\cT$ is represented by $\cT'=[t_{i',j',k'}]=\cT(P,Q,R)$.
 So
 \begin{equation}\label{changebas}
 \cT(P,Q,R)= P\otimes Q\otimes R(\cT):=\cT'=[t'_{i',j',k'}],\;
 t_{i',j',k'}=\sum_{i=j=k=1}^{m,n,l} p_{i'i}q_{j'j}r_{k'k}t_{i,j,k}.
 \end{equation}
 Clearly $\rank \cT(P,Q,R)=\rank \cT$.
 The following lemma is derived straightforward.
 \begin{lemma}\label{PQRchang}  Let $\cT\in\C^{m\times n\times l}, P\in\gl(m,\C),Q\in\gl(n,\C),
 R\in\gl(l,\C)$.  Let $\cT(P,Q,R)$ be defined as in (\ref{changebas}).  Then
 \begin{eqnarray*}
 &&\T_1(\cT(P,Q,R))=Q\T_1(\cT)R\trans, \quad \T_2(\cT(P,Q,R))=P\T_2(\cT)R\trans,\\
 &&\T_3(\cT(P,Q,R))=P\T_3(\cT)Q\trans.
 \end{eqnarray*}
 \end{lemma}

 For a finite dimensional space $\W$ of dimension $N$ denote by
 $\Gr(p,\W)$ the Grassmannian variety of all subspaces of dimension $1\le p\le N$.
 Recall that $\dim \Gr(p,\W)=p(N-p)$.  Let $\Gamma(p,\C^{m\times n})\subseteq \Gr(p,\C^{m\times n})$
 and $\Gamma(p,\rS(m,\C))\subseteq \Gr(p,\rS(m,\C))$ be the varieties
 of all $p$-dimensional subspaces in $\C^{m\times n}$ and $\rS(m,\C)$ that can be obtained as limit
 of $p$-dimensional subspaces in $\C^{m\times n}$ and $\rS(m,\C)$ respectively, which are
 spanned by rank one matrices.
 For integers $i\le j$ denote by $[i,j]$ the set of all integers $k, i\le k \le j$.
  The following result is known, e.g. \cite[Prop, 3.1, (iv)]{BL10}.
 \begin{lemma}\label{factsGp}  For $1<m,n\in\N$
 \begin{eqnarray}\label{factsGp1}
 \Gamma(p,\C^{m\times n})=\Gr(p,\C^{m\times n}) \textrm{ for } p\in [(m-1)(n-1)+1,mn],\\
 \label{factsGp2}
 \Gamma(p,\rS(m,\C))=\Gr(p,\rS(m,\C)) \textrm{ for } p\in [{m \choose 2}+1,{m+1 \choose 2}]
 \end{eqnarray}
 \end{lemma}
 \proof  To prove (\ref{factsGp1}) it is enough to show the case $p=(m-1)(n-1)+1$.
 Clearly $\Gamma(p,\C^{m\times n})$ is an irreducible variety of $\Gr(p,\C^{m\times n})$.
 It is left to show that $\dim \Gamma(p,\C^{m\times n})=\dim\Gr(p,\C^{m\times n})$.
 Let $\P V(r,m,n)\subset \P^{mn-1}$ be the projectivized variety of all matrices in $\C^{m\times n}\setminus\{0\}$
 of rank $r$ at most.  It is well known that $\dim \P V(r,m,n)= r(m+n-r)-1$, e.g. \cite{FK07}.
 Hence any generic projective linear subspace of dimension $(m-1)(n-1)$ in $\P^{mn-1}$ will intersect the Segre
 variety $\P V(1,m,n)= \P^{m-1}\times\P^{n-1}$ at a fixed number of points \cite[\S6]{FK07}
 \begin{equation}\label{degV1}
 \deg \P^{m-1}\times\P^{n-1}={m+n-2 \choose m-1}.
 \end{equation}
 Thus, a generic projective subspace spanned by
 $(m-1)(n-1)+1$ points on Segre variety, i.e. $(m-1)(n-1)+1$ rank one matrices.  Therefore
 \begin{eqnarray*}
 \dim  \Gamma((m-1)(n-1)+1,\C^{m\times n})= ((m-1)(n-1)+1)(m+n-2)=\\
 \dim \Gr((m-1)(n-1)+1,\C^{m\times n})\Rightarrow \Gamma(p,\C^{m\times n})=\Gr((m-1)(n-1)+1,\C^{m\times n}).
 \end{eqnarray*}

 To prove (\ref{factsGp2}) is enough to show the case $p={m\choose 2}+1$.
 Let $\P \rS(r,m)$ be the projectivized variety of all $A\in\rS(m,\C)\setminus\{0\}$ of rank $r$
 at most.  It is well known that $\codim \P \rS(r,m) = {m-r +1\choose 2}$ in $\P\rS(m,\C)$,
 e.g. \cite{FK07}.  Hence a generic projective linear subspace of dimension $m \choose 2$ in
 $\P\rS(m,\C)$ intersects $\P\rS(1,m)$ at a fixed number of points
 \cite{FK07}:
 \begin{equation}\label{degSym1}
 \deg\P\rS(1,m)=\prod_{j=0}^{m-2} {{m+j
 \choose m-1-j}\over {2j+1 \choose j}}.
 \end{equation}
 Similar arguments for the previous case show that $\Gamma({m \choose 2}+1,\rS(m,\C))$ and
 $\Gr({m \choose 2}+1,\rS(m,\C))$ have the same dimensions.  Hence the two varieties are equal.  \qed

 \begin{lemma}\label{eqdefbrank}  Let $0\ne \cT\in \C^{m\times n\times l}$.  Then $r\ge\brank \cT$,
 if there exists $\U\in\Gamma(r,\C^{m\times n})$ such that $\U\supseteq\T_3(\cT)$.  Furthermore, $r=\brank \cT$
 if there is no $\V\in\Gamma(r-1,\C^{m\times n})$ such that $\V\supseteq\T_3(\cT)$.  In particular,
 $\brank \cT\ge \dim\T_3(\cT)$.
 \end{lemma}
 \proof  Suppose that $\rank \cT=r$.  Then Theorem \ref{bascharrnk} yields the existence
 $\U\in\Gamma(r,\C^{m\times n})$ such that $\U\supseteq\T_3(\cT)$.  Suppose now that $\cT$ is a limit
 of $\cT'_q,q\in\N$ of rank $r'\le r$.  It is trivial to see that $\cT$ is a limit
 of $\cT_q,q\in\N$ of rank $r$.
 Then $\T_3(\cT_q)\subseteq U_q\in\Gamma(r,\C^{m\times n})$ for each $q\in\N$.
 Take a convergent subsequence $\U_{q_s}\to\U\in\Gamma(r,\C^{m\times n})$.  So $T_{k,3}(\cT)\in \U, k=1,\ldots,l$.
 Hence $\T_3(\cT)\subseteq \U$.  If $\brank\cT=r$ there is no $\V\in\Gamma(r-1,\C^{m\times n})$ such that
 $\V\supseteq\T_3(\cT)$.  Clearly, if $\U\supseteq \T_3(\cT)$ then $\dim \U\ge \dim \T_3(\cT)$.
 Hence $\brank\cT\ge \dim\T_3(\cT)$. \qed

 We now recall some basic results for matrices we need here.  Consult for example with \cite{MaM64,BFP95}.
 For $m\in\N$ let $\an{m}:=\{1,\ldots,m\}$.  For $k\in\an{m}$ denote by $2^{\an{m}}_k$ the set of subsets
 of $\an{m}$ of cardinality $k$.  Then $\alpha\in 2^{\an{m}}_k$ is viewed as
 $\alpha=\{\alpha_1,\ldots,\alpha_k\}$, where $1\le \alpha_1<\ldots<\alpha_k\le m$.
 Denote $ \|\alpha\|:=\sum_{j=1}^k \alpha_j,\alpha^c:=\an{m}\backslash\alpha$.
 For $A=[a_{i,j}]_{i=j=1}^{m,n}\in\C^{m\times n}, \alpha\in 2^{\an{m}}_k, \beta\in 2^{\an{n}}_l$
 denote $A[\alpha,\beta]:=[a_{\alpha_i,\beta_j}]_{i=j=1}^{k,l}\in \C^{k\times l}$.
 Recall that for $p\in\an{\min(m,n)}$ the $p-th$ compound of $A$, denoted as
 $\rC_p(A)\in \C^{{m \choose p}\times {n \choose p}}$, is a matrix whose rows and columns
 are indexed by $\alpha\in 2^{\an{m}}_p, \beta\in 2^{\an{n}}_p$ and its $(\alpha,\beta)$ entry
 is given by $\det A[\alpha,\beta]$.  If we view $A$ as a linear transformation $\hat A:\C^{n}\to
 \C^{m}$ given by $\x\mapsto A\x$, then $\rC_p(A)$ represents $\bigwedge^p \hat A:
 \bigwedge^p \C^n \to \bigwedge^p\C^m$ in the corresponding bases.
 Clearly
 $$\rC_p(A\trans)=\rC_p(A)\trans, \quad \rC_p(I_m)=I_{m\choose p}.$$
 (Here $I_m$ is $m\times m$ identity matrix.)
 The Cauchy-Binet formula yields
 \begin{equation}\label{CBform}
 \rC_p(AB)=\rC_p(A)\rC_p(B) \textrm{ for any } A\in \C^{m\times n}, B\in\C^{n\times l}.
 \end{equation}
 Let $m=n$.
 For $p\in\an{m-1}$ denote by $\rC_{-p}(A)$ a matrix whose rows and columns
 are indexed by $\alpha\in 2^{\an{m}}_p, \beta\in 2^{\an{m}}_p$ and its $(\alpha,\beta)$ entry
 is given by $(-1)^{\|\alpha\|+\|\beta\|}\det A[\alpha^c,\beta^c]$.  So $\rC_{-1}(A)\trans$ is the adjoint
 of $A$, denoted as $\adj A$.  The Laplace expansion yields that
 \begin{equation}\label{Cpident}
 \rC_p(A)\rC_{-p}(A)\trans=\rC_{-p}(A)\trans\rC_p(A)=(\det A)\; I_{m \choose p}, \quad p=1,\ldots,m-1.
 \end{equation}
 Hence
 \begin{equation}\label{C-pident}
 \rC_p(A^{-1})=\frac{1}{\det A} \rC_{-p}(A)\trans \textrm{ for } A\in\gl(m,\C).
 \end{equation}
 We now recall Strassen's result \cite{Str}.
 \begin{theo}\label{strass} Let $\cT=[t_{i,j,k}]_{i=j=k}^{3}\in\C^{3\times 3\times 3}$.
 Denote $X_k:=[t_{i,j,k}]_{i=j=1}^3\in\C^{3\times 3}, k=1,2,3$.  Let $f(X_1,X_2,X_3):=\det (X_1(\adj X_2)X_3
 -X_3(\adj X_2) X_1)$ be a polynomial of degree $12$ in the entries of the matrices $X_1,X_2,X_3$.
 Then $f(X_1,X_2,X_3)=s(X_1,X_2,X_3)\det X_2$.  The variety of all $\cT$ of border rank $4$ at most
 is a hypersurface in $\C^{3\times 3\times 3}$ of degree $9$ given by equation $s(X_1,X_2,X_3)=0$.
 \end{theo}

 The following result is straightforward.
 \begin{lemma}\label{solhomeq}.  Let $A\in\C^{m\times n}$ and assume that $\rank A\le k<n$.
 Fix $\alpha\in 2^{\an{m}}_{k},\beta=\{\beta_1,\ldots,\beta_{k+1}\}\in 2^{\an{n}}_{k+1}$.
 Let $\x(\alpha,\beta)=(x_1,\ldots,x_n)\trans\in\C^n$ be defined as follows.  $x_j=0$ if $j\not\in\beta$.
 If $j=\beta_i$ then $x_j=(-1)^{i-1}\det A[\alpha,\beta\backslash\{\beta_i\}]$.  Then $A\x(\alpha,\beta)=\0$.
 Furthermore, $\x(\alpha,\beta)\ne \0$ for some $\alpha,\beta$ if and only if $\rank A=k$.
 \end{lemma}

 The following result is known \cite{LM06}.
 \begin{theo}\label{commutcond}  Let $\cT\in\C^{m\times m\times l}$.
 Assume that $\brank\cT\le m$.  Then for any $A,B,C\in\T_3(\cT)$ the following equalities
 hold.
 \begin{equation}\label{commutcond1}
 \rC_p(A)\rC_{-p}(B)\trans\rC_p(C)=\rC_p(C)\rC_{-p}(B)\trans\rC_p(A) \textrm{ for } p=1,\ldots,m-1.
 \end{equation}
 \end{theo}
 \proof Assume that $\rank\cT=m$ and $\T_3(\cT)$ contains an invertible matrix $B$.
 So $\U=\span(\u_1\v_1\trans,\ldots,\u_m\v_m\trans)\supset\T_3(\cT)$.
 Hence $B=\sum_{i=1}^m b_i\u_i\v_i\trans$.  Since $\rank B=m$ we have that $b_i\ne 0, i=1,\ldots,m$,
 and $\u_1,\ldots,\u_m$ and $\v_1,\ldots,\v_m$ linearly independent.  There exist $P,Q\in\gl(m,\C)$
 such that $P\u_i=Q\v_i=\e_i:=(\delta_{i1},\ldots,\delta_{im})\trans, i=1,\ldots,m$.
 So any matrix in $P\U Q\trans$ is a diagonal matrix.  Hence
 $PAQ(PBQ)^{-1} PCQ=PCQ(PBQ)^{-1} PAQ$, i.e.
 \begin{equation}\label{commutcond2}
 AB^{-1}C=CB^{-1}A,
 \end{equation}
 for any $A,B,C\in\U$.
 Take the $p-th$ compound of the above equality, use (\ref{CBform}) and (\ref{C-pident})
 to deduce (\ref{commutcond1}) for any $A,B,C\in \U$.  Since almost any $B\in \U$ is invertible
 (\ref{commutcond1}) holds for any $A,B,C\in\U$.  In particular it holds for any $A,B,C\in\T_3(\cT)$.

 Assume now that $\brank \cT\le m$.  Then $\cT$ is a limit of $\cT_q,q\in\N$ of rank $m$.
 Furthermore, it is easy to see that we can assume that each $\T_3(\cT_q)$ contains an invertible matrix.
 Let $\T_3(\cT_q)\subset\U_q$, where $\U_q$ is a span of $m$ matrices of rank one.
 Hence any three matrices $A,B,C\in\U_q$ satisfy (\ref{commutcond1}).  Assume that $\U_q,q\in\N$ converges
 to $\U\in \Gamma(m,\C^{m\times n})$.  Then any $3$ matrices in $\U$ satisfy (\ref{commutcond1}).
 The proof of Lemma \ref{eqdefbrank} yields that $\T_3(\cT)\subset \U$.  Hence any $3$ matrices in
 $\T_3(\cT)$ satisfy (\ref{commutcond1}).  \qed

 The following result is well known, e.g. \cite{Fri08}.
 \begin{equation}\label{grank2mn}
 \brank\cT\le \min(n,2m) \textrm{ for any } \cT\in\C^{2\times m\times n} \textrm{ where } 2\le m\le n.
 \end{equation}
 \section{Some subspaces of singular matrices satisfying \eqref{XYZcomutcon}}
 For a subspace $\U\subset \C^{m\times n}$ define $\maxrank \U=\{\max \rank A, A\in\U\}$.
 The following theorem analyzes the condition (\ref{XYZcomutcon}) for a subspace $\U\subset\C^{m\times m}$
 satisfying $\maxrank \U=m-1$.
 \begin{theo}\label{mrank=m-1}  Let $\U\subset\C^{m\times m}$ and assume that $\maxrank \U=m-1$.
 Then any three matrices in $\U$
 satisfy (\ref{XYZcomutcon}) if and only if one of the following mutually exclusive conditions hold.
 \begin{enumerate}
 \item\label{mrank=m-1a}  There exists a nonzero $\u\in\C^m$ such that either $\U\u=\0$ or $\u\trans \U=\0\trans$.
 \item\label{mrank=m-1char} $m\ge 3,\dim \U=k+1\ge 2$.  There exists $P,Q\in\gl(m,\C)$ such that $P\U Q$ has a following basis $F_0,\ldots,F_k$.  The last row and column of $F_0,\ldots,F_{k-1}$ is zero, i.e. $F_i=G_i\oplus 0, F_i\in\C^{(m-1)\times (m-1)}, i=0,\ldots,k-1$,
 $G_0=I_{m-1}$, and
 \begin{equation}\label{formG1k}
 F_k=\left[\begin{array}{cc} G_k&\f\\ \g\trans&0\end{array}\right],\;G_k\in \C^{(m-1)\times (m-1)},
 \;\0\ne \f,\g\in\C^{m-1},\;\g\trans\f=0.
 \end{equation}
 Furthermore there exists two subspace $\X,\Y\subset\C^{(m-1)}$ with the following properties
 \begin{eqnarray}\label{condXY1}
 \f\in\X,\;\g\in\Y,\; \g\trans \X=\f\trans \Y=\0\trans,\;  G_k\X\subseteq \X,\;
 G_k\trans \Y\subseteq \Y,\\
 \label{condXY2}
 G_i\X=\0,\; G_i\trans \Y=\0,\; i=1,\ldots,k-1.
 \end{eqnarray}
 \end{enumerate}
 \end{theo}
 \proof  Let $A\in\U, \rank A=m-1$.  Then $\adj A=\u(A)\v(A)\trans$ for some nonzero $\u(A),\v(A)\in\C^m$.
 Since $A (\adj A)=(\adj A)A=0$ we deduce that $A\u(A)=A\trans \v(A)=\0$.  Suppose first that $\U\u=\0$
 for some nonzero $\u\in\C^m$.  So for each $A\in\U, \rank A=m-1$ we must have that $\span(\u(A))=\span(\u)$.  So we may assume that $\u(A)=\u$.  Hence for any $B\in \U$ $B\adj(A)=0$.
 Since $\adj(Y)=0$ if $\rank Y<m-1$, we deduce that any three matrices in $\U$ satisfy (\ref{XYZcomutcon}). Similarly, (\ref{XYZcomutcon}) holds if there exists nonzero $\u$ such that $\u\trans\U=\0$.

 Assume now that condition \ref{mrank=m-1a} does not hold.
 Then for most of $B\in\U$
 \begin{equation}\label{Buvne0}
 B\u(A)\ne \0, \quad B\trans\v(A)\ne \0.
 \end{equation}
 Assume now that (\ref{XYZcomutcon}) holds.
 Let $X=B,Y=A,Z=C$, and $B,C$ satisfy (\ref{Buvne0}).  Then
 \[\span(B\u(A))=\span(C\u(A))=\span(\x(A)), \span(B\trans \v(A))=\span(C\trans \v(A))=\span(\y(A)),\]
 for some nonzero $\x(A),\y(A)\in\C^m$.
 Hence, there exists two nontrivial linear functionals $\phi,\psi:\U\to \C$, depending on $A$, such that
 \[B\u(A)=\phi(B)\x(A),\quad B\trans \v(A)=\psi(B)\y(A) \textrm{ for all } B\in\U.\]
 Using (\ref{XYZcomutcon}) for $X=B,Y=A,Z=C$ and the above assumptions we obtain the equality
 $\phi(B)\psi(C)=\phi(C)\psi(B)$.  Choosing $B,C$ satisfying (\ref{Buvne0}) we get that
 $\frac{\phi(B)}{\psi(B)}=\frac{\phi(C)}{\psi(C)}$.  Hence $\psi(B)=a\phi(B)$ for all $B\in\U$, and
 $a\ne 0$.  By replacing $\y(A)$ by $a\y(A)$ we may assume that $\psi=\phi$.  Hence
 for each $A\in\U, \rank A=m-1$ we have the equality
 \begin{equation}\label{mrank=m-1c}
 B\u(A)=\phi_A(B)\x(A),\;B\trans \v(A)=\phi_A(B)\y(A) \textrm{ for all } B\in\U,
 \end{equation}
 for a corresponding nontrivial linear functional $\phi_A:\U\to\C$.

 Choose $P,Q\in \gl(m,\C)$ such that $F_0=I_{m-1}\oplus 0\in \U'=P\U Q$.
 Note that $\u(F_0)=\v(F_0)=\e_m=(\delta_{1m},\ldots,\delta_{mm})\trans$.  Then there exists a nonzero linear functional $\phi_0:\U'\to \C$
 such that $B\e_m=\phi_0(B)\x_0, B\trans \e_m=\phi_0(B)\y_0$ for some nonzero $\x_0,\y_0\in\C^m$ and $B\in\U'$.  Observe that for any $B\in\C^{m\times m}$ we have $\det(F_0+tB)=tb+O(t^2)$, where $b$ is the
 $(m,m)$ entry of $B$.  Hence $b=0$ for $B\in \U'$.
 Since $(m,m)$ entry of $B$ is zero it follows that $\x_0,\y_0$ have the last coordinate zero, i.e.
 $\x_0\trans=(\f\trans,0),\y_0\trans=(\g\trans,0)$.
 Consider next the strict subspace of $\U_0\subset \U'$ satisfying $\phi_0(B)=0$.
 For $B\in\U_0$ we have that $B\e_m=B\trans\e_m=\0$, i.e. the last row and column matrices in $\U_0$
 are zero.  Clearly $F_0\in\U_0$.  Let $F_0,\ldots,F_{k-1}$ be a basis in $\U_0$.  Let $\phi_0(F_k)=1$.
 The assumption that $\phi_0(F_k)=1$ yields that $F_k$ is of the form given in (\ref{formG1k}).
 (We will show the condition $\g\trans\f=0$ later.)
 So $F_0,\ldots,F_k$ is a basis of $\U'$.
 Let
 \[F(\z)=F_0+\sum_{i=1}^{k} z_iF_i,\quad G(\z)=I_{m-1}+\sum_{i=1}^k z_i G_i, \quad \z=(z_1,\ldots,z_k)\trans\in\C^k.\]
 Assume that $G(\z)$ is invertible.
 A straightforward calculation shows
 \begin{eqnarray}\notag
 &&\u(F(\z))\trans=(z_k\f(\z)\trans, -1),\; \v(F(\z))\trans=(\det G(\z))(z_k\g(\z)\trans,-1),\\
 &&\textrm{where } \f(\z)=G(\z)^{-1}\f,\;\g(\z)\trans=\g\trans G(\z)^{-1},
 \textrm{ and } \g\trans G(\z)^{-1}\f=0.\label{explforuvz}
 \end{eqnarray}
 Indeed, the equalities $F(\z)\u(\z)=F(\z)\trans\v(\z)=\0$ are verified straightforward.
 The equality $\g\trans G(\z)^{-1}\f=0$ must hold if $z_k\ne 0$.  The continuity argument yields this
 condition for $z_k=0$ if $G(\z)$ is invertible.   Note that the condition $\g\trans G(\z)^{-1}\f=0$ for $\z=\0$ yields that $\g\trans\f=0$. To see that $\adj F(\z)=\u(\z)\v(\z)\trans$ just observe that the $(m,m)$ entry of $\adj F(\z)=\det G(\z)$.
 Observe next that (\ref{mrank=m-1c}) yields
 \begin{equation}\label{idfrmlem2.8}
 G(\w)\u(\z)=\phi_\z(\w)\x(\z),
 \end{equation}
 for some nonzero affine functional $\phi_\z(\w)$.  ($\phi_\z(\w)$ is affine since in the definition of $F(\z)$ the coefficient of $A_0$ is $1$.)
 If $z_k\ne 0$ we claim that we can choose
 \begin{equation}\label{choicexz}
 \x(\z)\trans=(\f(\z)\trans,0)
 \end{equation}
 Indeed, chose $\w=\0$ so $G(\0)=I_{m-1}$.  Clearly $F(\0)\u(\z)=(z_k\f(\z)\trans,0)\trans=z_k(\f(\z)\trans,0)\trans$.
 Let $\hat \w=(w_1,\ldots,w_{k-1},0)\trans$.  Then we get the equality
 \begin{equation}\label{Fwxzid}
 F(\hat \w)\u(\z)=\left[\begin{array}{c}z_k G(\hat\w)\f(\z)\\0\end{array}\right]=
 \phi_z(\hat \w)\left[\begin{array}{c}\f(\z)\\0\end{array}\right].
 \end{equation}
 If $z_k\ne 0$ then $\f(\z)$ is an eigenvector of $G(\hat\w)$ for each $\hat\w$.
 Assume that $G(\z)$ is invertible for any $\z$ satisfying
 satisfying $\|\z\|_{\max}\le r$ for some $r>0$.
 Use the continuity argument to deduce that $\f(\z)$ is an eigenvector of $G(\hat\w)$ for any
 $\z$ satisfying satisfying $\|\z\|_{\max}\le r$.  Letting $\z=\0$ we get that $\f=\f(\0)$ is an eigenvector for each $G(\hat \w)$.  Hence
 $G_i\f=\lambda_i\f$ for $i=0,\ldots,k-1$, where $\lambda_0=1$.
 By replacing $G_i$ with $G_i-\lambda_iI_{m-1}$ we may assume without loss of generality that
 $G_i\f=0$ for $i=1,\ldots,k-1$.  Let $\|\z\|_{\max}\le r$.  Since $f(\z)$ is an eigenvector
 of $G_i$ and $\f(0)$ corresponds to the zero eigenvalue of $G_i$ it follows $G_i\f(\z)=\0$, were $i=1,\ldots,k-1$.  Let $\z=(0,\ldots,0,z_k)\trans$ and $|z_k|< r$.
 So
 \[\f(\z)=(I_{m-1}+z_kG_k)^{-1}\f=\sum_{j=0}^{\infty}(-z_k)^jG_k^j\f.\]
 Let $\X,\Y$ be the cyclic subspaces spanned by $G_k^j\f,j=0,\ldots,$ and $(G_k\trans)^j\g,j=0,\ldots,$
 respectively.  Clearly, $G_k\X\subseteq \X, G_k\trans \Y\subseteq\Y$.
 The condition that $G_i\f(\z)=\0$ yields that $G_i\X=\0$ for $i=1,\ldots,k-1$.
 So $G(\hat\w)\f(\z)=\f(\z)$.
 The condition $\g\trans \f(\z)=0$ yields that $\g\trans \X=\f\trans\Y=\0\trans$.
 Observe next that (\ref{Fwxzid}) yields that $\phi_\z(\hat\w)=z_k$.  In view of
 (\ref{mrank=m-1c}) it follows that
 \[F(\hat\w)\trans \v(\z)=z_k\y(\z), \quad \y(\z)\trans=a(\z)(\g(\z)\trans,0) \textrm{ for some } 0\ne a(\z)\in\C.\]
 Hence $G_i\trans \g(\z)=0$ and $G_i\trans\Y=\0$ for $i=1,\ldots,k-1$.
 This establishes the conditions \emph{2} of the theorem.

 Vice versa, suppose that  the conditions \emph{2} of the theorem hold.
 Let $\U'=P\U Q$.
 Define $F(\z),G(\z),\u(\z),\v(\z),\f(\z),\g(\z)$ as above.
 It is enough to show the condition (\ref{mrank=m-1c}),
 where $A=F(\z), B=F(\w), C=F(\w'))$ and $\det G(\z)\ne 0$. Observe next that
 (\ref{condXY1}-\ref{condXY2}) yield that
 \[f(\z)=(I_{m-1}+z_kG_k)^{-1}\f,\; \g(\z)\trans =\g\trans (I_{m-1}+z_kG_k)^{-1}, \g\trans \f(\z)=\f\trans \g(\z)=0.\]
 Then
 \begin{eqnarray*}
 F(\w)\u(\z)=((z_k(I_{m-1}+w_k G_k)(I_{m-1}+z_k G_k)^{-1}\f-w_k\f)\trans,0)\trans=\\
 (z_k-w_k)(\f(\z)\trans,0)\trans.
 \end{eqnarray*}
 Similarly,
 \[\v(\z)\trans F(\w')=(z_k-w'_k)\det G(\z)(\g(\z)\trans,0).\]
 Hence the condition (\ref{XYZcomutcon}) holds for $Y=F(\z),X=F(\w), Z=F(\w')$
 when $\det G(\z)\ne 0$.  Since most of $F(\z)\in \U'$ satisfy the condition that $\det G(\z)\ne 0$
 we deduce that (\ref{XYZcomutcon}) for each $X,Y,Z\in\U'$.
 \begin{theo}\label{charrnk3ten}  A tensor $\cT\in\C^{3\times 3\times 3}$ has border rank $3$ at most
 if and only if $T_p(\cT)$ satisfies the condition (\ref{XYZcomutcon}) for some $p\in\{1,2,3\}$.
 \end{theo}
 \proof  Theorem \ref{commutcond} implies that each $\T_p(\cT)$ satisfies the condition (\ref{XYZcomutcon}).
 It is enough to consider the case where $\T_3(\cT)\ne\{0\}$ satisfies the condition (\ref{XYZcomutcon}).
 Suppose first that $\dim\T_p(\cT)\le 2$ for some $p\in\{1,2,3\}$.  Then by changing basis in the $p$-th component of $\C^3\otimes\C^3\otimes \C^3$ and interchanging the first and the $p$-th component, we can assume that $\cT$ as $2\times 3\times 3$ tensor.  \eqref{grank2mn} yields that $\brank\cT\le 3$.

 Assume now that $\dim\T_p(\cT)=3$ for $p=1,2,3$.
 Suppose first that $\max\rank\T_3(\cT)=1$.  So $\T_3(\cT)$ has a basis consisting of rank one matrices.  Theorem \ref{bascharrnk} implies that $\rank\cT=3$, hence $3\ge\brank\cT$.

 Assume now that $\maxrank\T_3(\cT)=3$, i.e. there exists an invertible $Y\in\T_3(\cT)$.  By considering $P=Y^{-1}$ and changing a basis
 in the first factor of $\C^3\otimes\C^3\otimes\C^3$ we may that $Y=I\in \T_3(\cT)$.
 Let $\T_3=\span(I,A_1,A_2)$.  So $A_1A_2=A_2A_1$.
 Recall that the variety of all commuting pairs $(X_1,X_2)\in (\C^{3 \times 3})^2$
 is irreducible \cite{MT}.  Hence a pair $(A_1,A_2)$ is a limit of
 generic commuting pairs $(X_1,X_2)$.  For a generic pair, $X_1$ has $3$ distinct eigenvalues.
 So $X_2$ is a polynomial in $X_1$.
 Thus, there exists $Q\in \gl(3,\C)$ such that $Q^{-1}\span (X_1,X_2)Q$ is a two dimensional subspace
 of $3\times 3$ diagonal matrices $\D\subset\C^{3\times 3}$.  Clearly $I\in\D$ and $\D$ is spanned by
 $3$ rank one diagonal matrices.  Hence $\span (I,X_1,X_2)\subseteq Q^{-1}\D Q\in\Gamma(3,\C^{3\times 3})$.
 Thus $\T_3(\cT)\subseteq\U\in \Gamma(3,\C^{3\times 3})$, and $\brank\cT\le 3$.

 Assume now that $\maxrank\T_3(\cT)=2$.  We claim that there is no nonzero $\u\in\C^3$ such that
 either $\T_3(\cT)\u=\0$ or $\u\trans\T_3(\cT)=\0\trans$.  Assume to the contrary that
 $\u\trans\T_3(\cT)=\0\trans$ for some nonzero $\u$.  By change of basis in the first component
 of $\C^3$ we may assume that $\u=\e_3$.  Hence the third row of each matrix in $\T_3(\cT)$ is zero.
 Hence $T_{3,1}=0$, i.e. $\dim\T_1(\cT)\le 2$ contradicting our assumptions.
 Similarly there is no nonzero $\u$ such that $\T_3(\cT)\u=\0$.
 Hence $\U=\T_3(\cT)$ satisfies the condition of part \emph{2} of Theorem \ref{mrank=m-1}.
 Since $m=3$ it follows that the subspaces $X,Y\subset \C^3$ are one dimensional.  By changing bases
 in $\C^2$ we can assume that $\f=\e_1,\g=\e_2$.  Hence the three $3$-slices of $\cT$ are
 \[F_0=\left[\begin{array}{ccc}1&0&0\\0&1&0\\0&0&0\end{array}\right], \;
 F_1=\left[\begin{array}{ccc}0&1&0\\0&0&0\\0&0&0\end{array}\right],\;
 F_3=\left[\begin{array}{ccc}*&*&1\\0&*&0\\0&1&0\end{array}\right].\]
 Do the following elementary row and column operations on $F_2$ to bring it to the form
 $F_2'=\left[\begin{array}{ccc}0&0&1\\0&0&0\\0&1&0\end{array}\right]$.
 First subtract a mulitple of the third row from the second and first row.
 Then subtract a multiple of the third column from the first column.
 Apply the same row and column operations on $F_0$ and $F_1$ to obtains the three $3$-slices
 $F_0,F_1,F_2'$ of the tensor $\cT'$.  
 Consider the three $2$-slices $\cT'$:
 \[T_{1,2}=\left[\begin{array}{ccc}1&0&0\\0&0&0\\0&0&0\end{array}\right], \;
 T_{2,2}=\left[\begin{array}{ccc}0&1&0\\1&0&0\\0&0&1\end{array}\right],\;
 T_{3,2}=\left[\begin{array}{ccc}0&0&1\\0&0&0\\0&0&0\end{array}\right].\]
 Interchange the first two columns in each of the above matrices to obtain the
 matrices $A_1,I,A_2$.   Note that $A_1A_2=A_2A_1=0$.  The previous arguments show that
 $\brank\cT=\brank\cT'\le 3$.  \qed

 We conclude this section with the following proposition.
 \begin{prop}\label{dimT3=3}  Let $\cT\in \C^{m\times m\times m}$.  Assume that $\dim T_3(\cT)\le m-1$.
 Then any three matrices $A,B,C\in\T_k(\cT)$ satisfy the conditions (\ref{commutcond1}) for
 $p=1,m-1$ and $k=1,2$.  For $m=4$ the condition (\ref{commutcond1}) holds also for
 $p=2$ and $k=1,2$.
 \end{prop}
 \proof  By changing a basis in the last component of $\C^m\otimes\C^m\otimes\C^m$ we may assume that
 $T_{m,3}=0$.  Hence the last row of each $T_{i,1}$ and the last column of each $T_{i,2}$ is zero.
 Theorem \ref{mrank=m-1} yields that any three matrices in $\T_1(\cT),\T_2(\cT)$ satisfy the conditions
 (\ref{XYZcomutcon}).  Note that for any two matrices $A,B\in\T_k(\cT)$ either $A(\adj B)=0$ or $(\adj B)A=0$.
 Hence for any three matrices in $\T_k(\cT)$ (\ref{commutcond1}) holds for $p=m-1$ and $k=1,2$.

 Assume now that $m=4$.
 We now show that for any three matrices in $\cT_k(\cT)$ and $k=1,2$ we have the equality
 $\rC_2(A)\rC_{-2}(B)\trans \rC_2(C)=0$.  It is enough to show this identity for $k=1$.  Since the last row of $A\in\T_1(\cT)$ is zero, it follows that $\rC_2(A)$ has three zero rows labeled $(1,4),(2,4), (3,4)$.
 Hence the zero rows of $\rC_{-2}(B)$ are the rows $(1,2), (1,3), (2,3)$.  So $\rC_{-2}(B)\trans$ has
 three zero columns $(1,2),(1,3),(2,3)$.  A straightforward calculation shows that $\rC_2(A)\rC_{-2}(B)\trans \rC_2(C)=0$.
 Hence $\rC_2(A)\rC_{-2}(B)\trans \rC_2(C)=\rC_2(C)\rC_{-2}(B)\trans \rC_2(A)=0$. \qed

 For $m=4$ it seems to us that the condition \eqref{XYZcomutcon} always implies the conditions
 (\ref{commutcond1}) for $p=2,3$.

 \section{Tensors in $\C^{(m-1)\times(m-1)\times m}$ of border rank $m$}
 Let $\cT\in \C^{(m-1)\times(m-1)\times m}$ be of rank $m$.  So
 $$\cT=\sum_{i=1}^m \u_i\otimes\v_i\otimes\w_i,\quad \u_i,\v_i\in\C^{m-1},\w_i\in\C^m,i=1,\ldots,m.$$
 We call $\cT$ of rank $m$ \emph{generic} if any $m-1$ vectors out $\u_1,\ldots,\u_m$ and $\v_1,\ldots,\v_m$
 are linearly independent.
 \begin{lemma}\label{symmet}  Let $\cT=[t_{i,j,k}]\in\C^{(m-1)\times(m-1)\times m}$ be a generic rank
 $m$ tensor.  Then there exists unique $L,R\ne 0$ (up to a nonzero scalars), such that
 \begin{eqnarray}\label{symmet2}
 L\T_3(\cT)\subset\rS(m-1,\C),\; \T_3(\cT)R\subset\rS(m-1,\C),\\ L R\trans=R\trans L=
 (\frac{1}{m-1} \trace(LR\trans))I_{m-1}.\label{symmet1}
 \end{eqnarray}
 Furthermore $L,R\in\gl(m-1,\C)$.
 \end{lemma}
 \proof  Let $U,V\in\gl(m-1,\C)$ such that $U\u_i=V\v_i=\e_i, i=1,\ldots,m-1$.
 Let $U\u_m=\x=(x_1,\ldots,x_{m-1})\trans, V\v_m=\y=(y_1,\ldots,y_{m-1})\trans$.
 The assumption that any $m-1$ vectors from $\u_1,\ldots,\u_m$ and $\v_1,\ldots,\v_m$ are linearly
 independent imply that all the coordinates of $\x$ and $\y$ are nonzero.
 Hence there exists a diagonal $D\in\gl(m-1,\C)$ such that $D\x=D^{-1}\y$.
 So $(D\e_i)(D^{-1}\e_i)\trans, i=1,\ldots,m-1$ are $m-1$ commuting diagonal matrices.
 Furthermore the matrix $(D\x)(D^{-1}\y)\trans$ is symmetric.  Hence $DU\T_3(\cT)V\trans D^{-1}
 \subset\rS(m-1,\C)$.  Thus
 $L=V^{-1}D^2 U$ and $R=(L^{-1})^{\trans}$ will satisfy the conditions of the lemma.
 It is left to show that $L$ and $R$ are unique up to a  multiple of a nonzero constant.
 For that we may assume already that  $\T_3(\cT)$ is spanned by $\e_i\e_i\trans, i=1,\ldots,m-1$
 and $\z\z\trans$ for some $\z$ with nonzero coordinates.  The assumptions that $L\e_i\e_i\trans$
 is symmetric for $i=1,\ldots, m-1$ yields that $L$ is a diagonal matrix.
 The assumption that $L\z\z\trans$ is symmetric implies that $L=dI_{m-1}$.  So if $L\ne 0$
 then it is a nonzero multiple of $I_{m-1}$.  Similar results hold for $R$.  In particular, $R\trans$
 is an inverse of $L$ times a nonzero constant.  \qed
 \begin{lemma}\label{brankmt(m-1)(m-1)m}  Let $\cT=[t_{i,j,k}]\in\C^{(m-1)\times(m-1)\times m}$ be a border rank $m$ at most.  Then there exist $L,R\in\C^{(m-1)\times (m-1)}\backslash\{0\}$ such that (\ref{symmet1}) holds.
 \end{lemma}
 \proof  There exist a sequence of $\cT_k\in\C^{(m-1)\times (m-1)\times m}$ of rank $m$ at most that
 converge to $\cT$.  By perturbing each $\cT_k$ we can assume that each $\cT_k$ ia generic tensor
 of rank $m$.  So there exists $L_k,R_k\in\gl(m-1,\C)$ satisfying (\ref{symmet1}).
 Normalize $L_k,R_k$ to have $\trace (L_kL_k^*)=\trace(R_kR_k^*)=1$.  Since the set $\{A\in\C^{(m-1)\times (m-1)}, \trace (A A^*)=1\}$ is compact, there exists a subsequence $k_{p},p\in\N$,
 such that $L_{k_p}\to L, R_{k_p}\to R$ and $\T_3(\cT_{k_p})$ converges to $\U\in \Gamma(m,\C^{m\times m})$.
 Clearly $L\U,\U R\subset \rS(m,\C)$, and $L,R$ satisfy the equality in (\ref{symmet2})-(\ref{symmet1}).
 As $\U\supseteq \T_3(\cT)$ we deduce the lemma.  \qed
 \begin{lemma}\label{sufconsym}  Let $\cT=[t_{i,j,k}]\in\C^{(m-1)\times(m-1)\times r}$, where $3\le r$.
 Denote by $T_k:=[t_{i,j,k}]_{i=j=1}^{m-1}$ the $(k,3)$-slice of $\cT$ for $k=1,\ldots,r$.  Then the two systems
 \begin{eqnarray}\label{sufconsym1}
 T_kR-R\trans T_k\trans =0,\;k=1,\ldots,r, R\in\C^{(m-1)\times(m-1)},\\
 LT_k-T_k\trans L\trans=0,\;k=1,\ldots,r, L\in\C^{(m-1)\times(m-1)}\label{sufconsym2}
 \end{eqnarray}
 have nontrivial solutions $R,L$ if and only if the following conditions hold. Let $C_R(T_1,\ldots,T_r)$,
 $C_L(T_1,\ldots,T_r) \in \C^{r(m-1)^2\times (m-1)^2}$ be the coefficient matrices of the systems
 (\ref{sufconsym1}) and (\ref{sufconsym2})
 in $(m-1)^2$ variables, (the entries of $R$ and $L$ respectively), and $r{m-1\choose 2}$ equations.
 Then $\rank C_R(T_1,\ldots,T_r) < (m-1)^2$ and  $\rank C_L(T_1,\ldots,T_r) < (m-1)^2$.
 Equivalently any $(m-1)^2\times (m-1)^2$ minors of $C_R(T_1,\ldots,T_r)$ and  $C_L(T_1,\ldots,T_r)$
 vanishes.  This assumption is equivalent to the assumption that the entries of $T_1,\ldots,T_r$ satisfy corresponding system of $2{\frac{r(m-1)(m-2)}{2}\choose (m-1)^2}$ homogeneous polynomial equations of degree $(m-1)^2$.

 Assume furthermore that $\rank C_R(T_1,\ldots,T_r) =\rank C_L(T_1,\ldots,T_r) = (m-1)^2 -1$.
 Then nonzero solutions $R,L$ of (\ref{sufconsym1}) and (\ref{sufconsym2}) are unique up to multiples by
 nonzero constants.  The equalities (\ref{symmet1}) are equivalent to
 $2{\frac{r(m-1)(m-2)}{2}\choose (m-1)^2-1}^2(m-1)^2$ homogeneous polynomial equations of degree $2((m-1)^2-1)$.
 \end{lemma}
 \proof
 As $X-X\trans$ is a skew symmetric matrix, the condition that $X\in\rS(m-1,\C)$ is equivalent to the fact
 that the entries of $X$ satisfy $m-1\choose 2$ linearly independent conditions.
 So $C_L(T_1,\ldots,T_r),C_R(T_1,\ldots,T_r)\in\C^{r{m-1\choose 2}\times (m-1)^2}$.
 Note that any element of  $C_R(T_1,\ldots,T_r)$ and $C_L(T_1,\ldots,T_r)$ is a linear function
 in the entries of some matrix $T_k$.
 Hence any $(m-1)^2\times (m-1)^2$ minor of $C_R(T_1,\ldots,T_l)$
 and $C_L(T_1,\ldots,T_l)$ is a polynomial of degree $(m-1)^2$ in entries of $\cT$.
 There are ${\frac{r(m-1)(m-2)}{2}\choose (m-1)^2}$ distinct minors of order $(m-1)^2$ of
 $C_R(T_1,\ldots,T_r)$ and $C_L(T_1,\ldots,T_r)$ respectively, which corresponds to a choice
 of $(m-1)^2$ rows from $r{m-1 \choose 2}$ rows.  Hence the total number of polynomial conditions for the existence
 of nonzero solution of (\ref{sufconsym1}) and (\ref{sufconsym2}) is equivalent to the vanishing of
 all $2{\frac{r(m-1)(m-2)}{2}\choose (m-1)^2}$ minors of  $C_R(T_1,\ldots,T_r)$
 and $C_L(T_1,\ldots,T_r)$ of order $(m-1)^2$.

 Suppose now that $\rank C_R(T_1,\ldots,T_r) =\rank C_L(T_1,\ldots,T_r) = (m-1)^2 -1$.
 Choose a solution for $L$ and $R$ as in Lemma \ref{solhomeq}.  If either $L$ or $R$ are zero matrices
 then (\ref{symmet1}) holds trivially.
 If $R,L\ne 0$ then the conditions (\ref{symmet1}) are equivalent to $2(m-1)^2$
 polynomial identities of degree $2((m-1)^2-1)$ in the entries of $T_1,\ldots,T_r$.
 The number of choices of $L$ and $R$ as described in Lemma \ref{solhomeq}
 is ${\frac{r(m-1)(m-2)}{2}\choose (m-1)^2-1}^2$ respectively.
 \qed

 We now discuss in detail the cases $m=4$ and $r=3,4$.  The case $r=3$ is the Strassen condition.
 \begin{theo}\label{equicstrcon}  Let $\cT=[t_{i,j,k}]_{i=j=k}^3\in\C^{3\times 3\times 3}$.
 Denote by $T_1,T_2,T_3\in\C^{3\times 3}$ the three $3$-slices of $\cT$.
 Let $C_{R}(T_1,T_2,T_3), C_L(T_1,T_2,T_3)\in \C^{9\times 9}$ be the matrix coefficients of the systems
 (\ref{sufconsym1}) and  (\ref{sufconsym2}) in the $9$ entries of $R$ and $L$ respectively.
 Then the border rank of $\cT$ is $4$ at most if and only if one of the following condition hold.
 \begin{enumerate}
 \item\label{equicstrcon1} $\det C_R(T_1,T_2,T_3)=0$.
 \item\label{equicstrcon2} $\det C_L(T_1,T_2,T_3)=0$.
 \end{enumerate}
 Equivalently, for any $\cT\in \C^{3\times 3\times 3}$
 \begin{equation}\label{equicstrcon3}
 \det C_R(T_1,T_2,T_3)=as(T_1,T_2,T_3),\quad \det C_L(T_1,T_2,T_3)=bs(T_1,T_2,T_3)
 \end{equation}
 for some nonzero $a,b\in \C$,
 where $s(T_1,T_2,T_3)$ is the Strassen polynomial described in Theorem \ref{strass}.
 \end{theo}
 \proof  Suppose first that $T_1,T_2,T_3\in\rS(3,\C)$.
 Assume that $T_1,T_2,T_3$ are three generic matrices.  Add a generic matrix $T_4\in\rS(3,\C)$.
 The proof of Lemma \ref{factsGp} yields that $\span(T_1,T_2,T_3,T_4)$ is spanned by $4$ rank one
 symmetric matrices.  Theorem \ref{bascharrnk} yields that $\rank \cT\le 4$.  Assume that $T_1(\adj T_2) T_3\ne T_3(\adj T_2)T_1$.  Theorem \ref{commutcond} implies that $\brank \cT\ge 4$.  Hence $\rank\cT=4$.
 Since any $T_1,T_2,T_3\in\rS(3,\C)$ can be approximated by three symmetric matrices in general position
 we deduce that $\brank\cT\le 4$ if the three $3$-slices of $\cT$ are symmetric matrices.
 Thus if (\ref{sufconsym1}) has a solution $R\in\gl(3,C)$ then $\brank\cT\le 4$.

 We now show that there exists
 $T_1,T_2,T_3\in\C^{3\times 3}$ such that (\ref{sufconsym1}) has only the trivial solution $R=0$.
 Let $T_1=I$, $T_2$ a diagonal matrix with $3$ distinct eigenvalues and $T_3=[s_{ij}]_{i=j=1}^3$, were all $s_{ij}\ne 0$.   The first condition of (\ref{sufconsym1}) yields that $R\in\rS(3,\C)$.
 The second condition of (\ref{sufconsym1}) imply that $R$ commutes with $T_2$.  Hence $R$ is a diagonal matrix $\diag(r_1,r_2,r_3)$.  The third condition of (\ref{sufconsym1}) is the condition
 $r_is_{ij}=s_{ji}r_j, i,j=1,\ldots,3$.  So if $s_{12}=s_{21}, s_{13}=s_{32}$ and $s_{23}\ne s_{32}$ it follows that $r_1=r_2=r_3=0$.

 On the other hand if $T_3$ is also a symmetric matrix with nonzero entries, then (\ref{sufconsym1})
 implies that $R=rI_3$.  Hence the condition $\det C_R(T_1,T_2,T_3)=0$ yield in the generic case, i.e. $\det R\ne 0$, that $\brank\cT\le 4$.
 By Strassen's theorem the set of $\cT\in\C^{3\times 3\times 3}$ of border rank $4$ is a hypersurface
 given by the equation $s(T_1,T_2,T_3)=0$.  Hence $\det C_R(T_1,T_2,T_3)=aS(T_1,T_2,T_3)$ for some $a\ne 0$.
 Similar results apply to $C_L(T_1,T_2,T_3)$.   \qed

 Assume $\cT\in\C^{3\times 3\times 3}$ has border rank $4$ at most.   Then $\det C_R(T_1,T_2,T_3)=0$,
 and $\det C_L(T_1,T_2,T_3)=0$.
 Then the choice of $R$ and $L$ given by Lemma \ref{solhomeq}
 is a column of $\adj C_R$ and $\adj C_L$ respectively.  So the entries of $R$ and $L$
 are homogeneous polynomials of degree $8$ in the entries of $\cT$.  Assume the generic case $\det R\ne 0$.
 Then the arguments in the proof of Theorem \ref{equicstrcon} show that (\ref{symmet1}) hold.
 Note that since each entry of $R$ and $L$ are polynomials of degree $8$ in the entries of $\cT$.
 So (\ref{symmet1}) are $18$ polynomial equations of degree $16$.  Since the only condition for
 $\cT\in\C^{3\times 3\times 3}$ is the vanishing of the Strassen polynomial, we deduce that
 each polynomial equation of (\ref{symmet1}) is given by the Strassen polynomial times a homogeneous polynomial of degree $7$.  In conclusion, in this case, (\ref{symmet1}) do not give any additional
 restriction on $\cT$.

 We now discuss the case $m=r=4$.
 So $\cT=[t_{i,j,k}]\in\C^{3 \times 3\times 4}$.
 We have four $3$-slices $T_{k}=[t_{i,j,k}]_{i=j=1}^3\in\C^{3\times 3}, k=1,\ldots,4$.
 Let $R=[x_{ij}]_{i=j=1}^3,L=[y_{ij}]_{i=j=1}^3$ be $3\times 3$ matrices with unknown entries.
 Then (\ref{sufconsym1}) and  (\ref{sufconsym2}) are $12$ equations homogeneous equations in $9$ variables $x_{11},\ldots,x_{33}$ and $y_{11},\ldots,y_{33}$, which are given by the coefficient matrices $C_R(T_1,T_2,T_3,T_4), C_L(T_1,T_2,T_3,T_4)\in \C^{12\times 9}$ respectively.  The condition that there exists nonzero $R$ and $L$ satisfying (\ref{sufconsym1}) and (\ref{sufconsym2}) respectively, are equivalent to the conditions $\rank C_R(T_1,T_2,T_3,T_4)\le 8, \rank C_L(T_1,T_2,T_3,T_4)\le 8$.  So each $9\times 9$ minor of
 $C_R(T_1,T_2,T_3,T_4)$, $C_L(T_1,T_2,T_3,T_4)$ is zero.  The number of these conditions is $2{12 \choose 9}=440$
 polynomial equations of degree $9$.  Fix submatrices $A,B\in\C^{9\times 9}$ of $C_R(T_1,T_2,T_3,T_4), C_L(T_1,T_2,T_3,T_4)$ respectively.  Then each column of $\adj A, \adj B$ respectively, represents
 a solution $R,L$ of (\ref{sufconsym1}) and  (\ref{sufconsym2}) respectively.  If $\rank A<8$ then $R=0$.
 If $\rank A=8$, then one of the $9$ columns of $\adj A$ is nonzero.  Similar conditions hold for $B$.
 So the number of the above choices of $R$ and $L$ is $9\times 220=1980$ for each of them.  Hence the total number
 of the above choices of pairs $R,L$ is $1980^2$.  For each choice of $R,L$ we assume that $18$ conditions
 given by (\ref{symmet1}) hold.  (To be precise, since $\tr LR\trans =\tr R\trans L$ we need at most
 $17$ equations of  (\ref{symmet1}).)  It is not known the the author if the conditions that $\rank C_R(T_1,T_2,T_3,T_4)\le 8, \rank C_L(T_1,T_2,T_3,T_4)\le 8$ imply (\ref{symmet1}), as in the case of $\cT\in\C^{3\times 3\times 3}$.

 \begin{theo}\label{334charbr4}$\cT=[t_{i,j,k}]_{i=j=k=1}^{3,3,4}\in\C^{3\times 3\times 4}$
 has a border rank $4$ at most if and only the following conditions hold.
 \begin{enumerate}
 \item\label{334charbr4a}
 Let $T_{k}:=[t_{i,j,k}]_{i=j=1}^3\in\C^{3\times 3},k=1,\ldots,4$ be the four $3$-slices of $\cT$.
 Then the ranks of $C_L(T_{1},\ldots,T_{4}),
 C_R(T_{1},\ldots,T_{4})$ are less than $9$.  (Those are $9-th$ degree equations.)
 \item\label{334charbr4b}
 Let $R, L$ be solutions of (\ref{sufconsym1})
 and (\ref{sufconsym2}) respectively as given in Lemma \ref{solhomeq}, (as described above).  Then (\ref{symmet1}) holds.
 (Those are $16-th$ degree equations.)
 \end{enumerate}
 \end{theo}
 \proof  Lemma \ref{brankmt(m-1)(m-1)m} implies that if $\brank\cT\le 4$ then the conditions
 \ref{334charbr4a}-\ref{334charbr4b} hold.  We now assume that the conditions
 \ref{334charbr4a}-\ref{334charbr4b} hold.
 Let $\U:=\T_3(\cT)$.
 Suppose first that $\dim \U\le 3$.  Pick $A_1,A_2,A_3\in \U$ such that
 $\span(A_1,A_2,A_3)=\U$.  Since each $A_i$ is a linear combination of $T_1,\ldots,T_4$, our assumption implies that there exists nonzero $R$ such that $A_iR- R\trans A_i\trans=0$ for $i=1,2,3$.  Hence $\det C_R(A_1,A_2,A_3)=0$ which is equivalent to the Strassen condition
 $s(A_1,A_2,A_3)=0$.  Strassen's theorem implies that $\rank \cT_3\le 4$.

 Assume now that $\dim\U=4$.  Lemma \ref{eqdefbrank} implies that $\brank\cT\ge 4$.  Let $R\in\C^{3\times 3}\backslash\{0\}$ be a solution of
 (\ref{sufconsym1}) for $m=4$.
 If $R\in\C^{3\times 3}$ has rank $3$ then $T_1':=T_1R,\ldots,T_4':=T_4R$ are $4$ linearly independent
 symmetric matrices.  Use Lemma \ref{factsGp} to deduce that
 $\cT'\in\C^{3\times 3\times 4}$ has border rank $4$.
 Similar results hold if $\rank L=3$.
 It is left to consider the case where $\max(\rank R,\rank L)\le 2$.
 We now consider a number of cases.\\

 \noindent
 A: $\rank C_R(T_1,\ldots,T_4)=\rank C_L(T_1,\ldots,T_4)=8$

 \noindent
 $\;$ I: $\rank L=\rank R=1$.  So after change of basis we can assume that $L=\e_3\e_3\trans$.
 Then the condition that $LT-T\trans L\trans  =0$ is equivalent to $T\trans \e_3=t\e_3$ for any $T\in\U$.
 We now consider the following  mutually exclusive subcases.

 \noindent
 $\;\;$ \emph{1}: $T_i\trans \e_3=0$ for $i=1,\ldots,4$.  Hence $\cT$ can be viewed as a tensor in
 $\C^{2\times 3 \times 4}$.  (\ref{grank2mn}) implies that $\brank\cT\le 4$.

 \noindent
 $\;\;$ \emph{2}: $\U$ contains $F_4:=\e_3\e_3\trans$.  So we can choose a basis $F_1,F_2,F_3,F_4$
 such that $F_i\trans\e_3=\0, i=1,2,3$.  Hence the tensor $\cT'\in \C^{3\times 3\times 3}$,
 whose three $3$-slices are $F_1,F_3,F_3$, can be viewed as a tensor in $\C^{2\times 3\times 3}$.
 (\ref{grank2mn}) implies that $\brank\cT'\le 3$ and
 and the border rank of $\cT$ is $4$ at most.

 \noindent
 $\;\;$ \emph{3}:  Let $T_k'$ be obtained from $T_k$ by deleting the last row for $k=1,2,3,4$.
 We claim that $T_1', \ldots,T_4'$ are linearly independent. Otherwise, there is a nontrivial combination
 $F\in\C^{3\times 3}$ of $T_1,\ldots,T_4$ such that the first two rows of $F$ are zero rows.
 Since $T_1,\ldots,T_4$ are linearly independent $F\ne 0$.  As $F\trans\e_3=t\e_3$ it follows that
 $F=t\e_3\e_e\trans, t\ne 0$.  This contradicts our assumption that the case \emph{2} does not hold.
 We now use the assumption that $TR-R\trans T\trans=0$
 and $R=\x\y\trans, \x=(x_1,x_2,x_3)\trans,\y=(y_1,y_2,y_3)\trans$ for each $T\in\U$.  So $T_k\x=s_k\y$ for $k=1,\ldots,4$.
 Suppose first that that
 all $s_k=0$.  Then we are done as in the case \emph{1}.  So we assume that $s_i\ne 0$ for some $i$.
 Since $R$ and $L$ have rank one, it follows that the condition (\ref{symmet1})
 implies that $R\trans L=LR\trans=0$.  Hence $x_3=\x\trans\e_3=0,y_3=\e_3\trans\y=0$.
 Let $\hat T_k\in\C^{2\times 2}$ obtained from $T_k$ by erasing the
 last row and column for $k=1,\ldots,4$.  Let $\hat\x\in\C^2$ be obtained from $\x\in\C^3$ by deleting the
 the last coordinate.  Then $\hat T_k\hat\x=s_k\hat \x$.  So by changing the coordinates in $\C^2$ we may assume
 that $\hat\x=(0,1)\trans$.  Combine the above conditions with the conditions that $T_i\trans \e_3=t_i\e_3, i=1,\ldots,4$  to deduce that there exists $P,Q\in\gl(3,\C)$ with the following properties.
 Let $\tilde T_k=PT_kQ=[\tilde t_{i,j,k}]_{i=j=1}^{3}\in \C^{3\times 3}, k=1,\ldots,4$.  Then
 \begin{equation}\label{ijzerocond}
 \tilde t_{i,j,k}=0 \textrm{ for }(i,j)=(1,2),\; (i,j)=(3,1),\; (i,j)=(3,2) \textrm{ and } k=1,\ldots,4.
 \end{equation}
 Take a generic subspace $\V\subset \C^{3\times 3}$ of dimension $4$ whose entries are zero at the places $(i,j)$
 given by (\ref{ijzerocond}).  We claim that $\V\in\Gamma(4,\C^{3\times 3})$.
 First take a matrix $D=[d_{ij}]_{i=j=1}^3\in\V$ such that $d_{ij}=0$ for $(1,1),(2,1),(2,2)$.
 Generically there would one matrix, up to multiplication by a scalar, such that $d_{33}\ne 0$.
 $D$ has rank one.  Now consider the $3$-dimensional subspace of $\V$ where the $(3,3)$ entry of each matrix
 is zero.  Then $\V$ can be viewed as a $3$-dimensional subspace in $\tilde \V\subset\C^{2\times 3}$.
 By Lemma \ref{factsGp} $\tilde\V\in \Gamma(3,\C^{2\times 3})$.  Hence
 $\V\in\Gamma(4,\C^{3\times 3})$  and $\brank \cT\le 4$.

 \noindent
 $\;$ II: $\max(\rank L,\rank R)=2$.
 By considering $\U\trans$ if necessary we may assume that $\rank L=2$.
 So there exist $P,Q\in \gl(3,\C)$ such
 that $PLQ=\diag(1,1,0)$.  Without loss of generality we may assume
 that $P=Q=I$.  Then each $LT_k$ is symmetric.  In particular
 $T_k \e_3=t_k\e_3$ and the $2\times 2$ submatrix $[t_{i,j,k}]_{i,j=1}^2$ is symmetric.
 We now claim that any four dimensional subspace $\V\subset \C^{3\times 3}$,
 such that each $T=[t_{ij}]_{i,j=1}^3\in\V$ satisfies $t_{12}=t_{21}, t_{13}=t_{23}=0$, is in $\Gamma(4,\C^{3\times 3})$.  As $\dim \V=4$ there exists $0\ne S=[s_{ij}]_{i=j=1}^3\in\V$ such that $0=s_{11}=s_{22}=s_{12}(=s_{21})$.
 Hence $\rank S=1$.
 For a generic $\V$ satisfying the above conditions $s_{33}\ne 0$.
 Consider now the $3$ -dimensional subspace $\W$ of $\V$ with $t_{33}=0$.
 Since $\W$ can be viewed as a $3$-dimensional subspace of $C^{3\times 2}$, Lemma \ref{factsGp}
 yields that $W\in \Gamma(3,\C^{3\times 3})$.  Hence $\V\in\Gamma(4,\C^{3\times 3})$
 and $\brank\cT\le 4$.\\

 \noindent
 B: $\min(\rank C_R(T_1,\ldots,T_4),\rank C_L(T_1,\ldots,T_4)<8$.
 By considering $\U\trans$ if necessary we can
 assume that $\rank C_L(T_1,\ldots,T_4)<8$.
 So there exist at least two linearly independent matrices $L_1,L_2\in\C^{3\times 3}$
 such that (\ref{sufconsym2}) holds.  If $\max(\rank L_1, \rank L_2)=3$ we deduce that $\brank\cT\le 4$ as in the beginning of our proof.  If $\max(\rank L_1, \rank L_2)=2$ we deduce that $\brank\cT\le 4$ as in the case A.II.
 So it is left to consider the
 case where $L_1$ and $L_2$ are rank one matrices such any their linear combination is also
 rank a one matrix.  It is easy to show that we can choose $P,Q\in \gl(3,\C)$ such that $PL_1Q=
 \e_3\e_3\trans$ and $PL_2Q$ is either $\e_2\e_3\trans$ or $\e_3\e_2\trans$.  So we have two cases.

 \noindent
 $\;$ I: $L_1=\e_3\e_3\trans, L_2=\e_2\e_3\trans$.  The condition (\ref{sufconsym2}) for $L_1$ yields
 that $T_k\trans \e_3=t_k\e_3$ for $k=1,2,3,4$,.  I.e. any $T\in\U$ has zero entries at the places
 $(3,1),(3,2)$. The condition (\ref{sufconsym2}) for $L_1$ yields $T_k\trans \e_3=t'_k\e_2$.
 Hence $t_k=t'_k=0$.  Thus the third row of each $T_k$ is zero.  So $\cT\in\C^{2\times 3\times 4}$
 and (\ref{grank2mn}) yields that $\brank\cT\le 4$.

 \noindent
 $\;$ II: $L_1=\e_3\e_3\trans, L_2=\e_3\e_2\trans$.  The condition (\ref{sufconsym2}) for $L_1$ yields
 that any $T\in\U$ has zero entries at the places $(3,1),(3,2)$.  The condition (\ref{sufconsym2})
 for $L_2$ yields  that $T_k\trans \e_2=t_k'\e_3$ for $k=1,2,3,4$.
 So the entries $(2,1),(2,2)$ are zero for each $T\in\U$.
 Take a nonzero $T_4'\in\U$ whose first row is zero.
 It is a rank one matrix.  Then either $(3,3)$ entry or $(2,2)$ entry of $T_4'$ is not equal to zero.
 Assume for simplicity of the argument that $(3,3)$ entry of $T_4'$ is nonzero.
 Hence $\U$ contains a three dimensional
 subspace  $\U'$ whose last row is zero.  Since a generic $3$ dimensional subspace of $2\times 3$
 matrices is spanned by rank one matrices it follows that $\cT$ has border rank $4$ at most.  \qed

 Note that in the proof of Theorem \ref{334charbr4} we used the condition \ref{334charbr4b},
 which are degree $16$ polynomial equations, only in the proof of the case A.I.3.
 Thus one can eliminate the use of degree $16$ polynomial equations, if one can show directly that
 a generic $4$-dimensional subspace of matrices satisfying (\ref{sufconsym1}) and (\ref{sufconsym2})
 for $R$ and $L$ of rank one, such that $RL\trans\ne 0$, is in $\Gamma(4,\C^{3\times 3})$.
 As an example, consider the case where $R=L=\e_3\e_3\trans$, which is ruled out by (\ref{symmet1}).
 Then the conditions (\ref{sufconsym1}) and (\ref{sufconsym2}) are equivalent to the assumption
 that $\T_3(\cT)$ is a four dimensional subspace of block diagonal matrices of the form
 \begin{equation}\label{blockdgcs}
 \left[\begin{array}{ccc}a&b&0\\c&d&0\\0&0&e\end{array}\right].
 \end{equation}
 Hence Theorem \ref{334charbr4} yields that $\T_3(\cT)\not\in\Gamma(4,\C^{3\times 3})$.
 It was shown in \cite{OB10} that the corresponding $\cT$ does not satisfy the degree $6$ polynomial
 equations found in \cite{LM04}.
 \section{Tensors in $\C^{4\times 4\times 4}$ of border rank 4 at most}
 \begin{theo}\label{charbr4}  $\cT=[t_{i,j,k}]_{i=j=k=1}^4\in\C^{4\times 4\times 4}$
  has a border rank $4$ at most if and only the following conditions hold.
 \begin{enumerate}
 \item\label{charbr4a} Any three matrices in $\T_1(\cT),\T_2(\cT),\T_3(\cT)$
 satisfy the conditions (\ref{XYZcomutcon}).
 (These are $5-th$ order degree equations on entries of $X,Y,Z$.)
 \item\label{charbr4b}
 For each $P_1,P_2,P_3\in \C^{4\times 4}$ let $\cT(P_1,P_2,P_3)=[t_{i,j,k}(P_1,P_2,P_3)]_{i=j=k=3}^4\in\C^{4\times 4\times 4}$
 be the tensor given by (\ref{changebas}).
 Let $S_{i_p,p}(P_1,P_2,P_3)\in \C^{4\times 4}, i_p=1,\ldots,4$  be the four $p$-slices
 of $\cT(P_1,P_2,P_3)$ for $p=1,2,3$.  (The entries of $S_{i_p,p}(P_1,P_2,P_3)\in\C^{4\times 4}$ are given
 by $t_{i_1,i_2,i_3}(P_1,P_2,P_3)$, where $i_p$ is fixed for a given $p\in\{1,2,2\}$ and $i_p\in\{1,2,3,4\}$.)
 Denote by $T_{i_p,p}(P_1,P_2,P_3)\in\C^{3\times 3}$ the submatrix obtained from $S_{i_p,p}(P_1,P_2,P_3)$
 by deleting the last row and column, for $i_p=1,2,3,4$.  Then
 \begin{eqnarray}\label{rankcondTP123a}
 \rank C_L(T_{1,p}(P_1,P_2,P_3),\ldots, T_{4,p}(P_1,P_2,P_3))\le 8,\\
 \rank C_R(T_{1,p}(P_1,P_2,P_3),\ldots,T_{4,p}(P_1,P_2,P_3))\le 8, \label{rankcondTP123b}
 \end{eqnarray}
 for $p=1,2,3$.  (Those are degree $9-th$ degree equations.)  Moreover the following conditions
 are satisfied for each $p\in\{1,2,3\}$.  Let $R_p(P_1,P_2,P_3), L_p(P_1,P_2,P_3)$ be solutions of (\ref{sufconsym1})
 and (\ref{sufconsym2}) respectively as given in Lemma \ref{solhomeq}.  Then (\ref{symmet1}) holds.
 (Those are degree $16-th$ degree equations.)
 \end{enumerate}
 \end{theo}

 To prove Theorem \ref{charbr4} we need to prove Corollary 5.6 of \cite{LM06}.
 \begin{theo}\label{landman}  Let $\cT\in\C^{4\times 4\times 4}$ and assume that any three matrices $X,Y,Z$ in $\T_p(\cT)$ satisfy (\ref{XYZcomutcon}) for $p=1,2,3$.  Then either $\brank \cT\le 4$ or
 or $\dim\T_p(\cT)\le 3, \dim\T_q(\cT)\le 3$ for two integers $1\le p<q\le 3$.
 Equivalently, by permuting factors in $\C^4\otimes\C^4\otimes \C^4$, if necessary, and changing bases in the first two components of $\C^4\times\C^4\times\C^4$ the tensor $\cT$ can be viewed as a tensor $\C^{3\times 3\times 4}$.
 \end{theo}
 The proof of this theorem is completed by considering a number of lemmas.
 \begin{lemma}\label{4invcase}  Let $\cT=[t_{i,j,k}]\in \C^{4 \times 4\times 4}$ and $p\in\{1,2,3\}$.
 Assume that $\T_p(\cT)$ contains an invertible matrix.  Then the condition
 (\ref{XYZcomutcon}) for any three matrices in $\T_p(\cT)$ implies that $\brank\cT\le 4$.
 \end{lemma}
 \proof  It is enough to consider the case $p=3$.
 Assume that  $Y\in\T_3(\cT)$ is invertible.  By considering $P=Y^{-1}$ and changing a basis
 in the first factor of $\C^4\otimes\C^4\otimes\C^4$ we may that $Y=I\in \T_3(\cT)$.
 Let $\T_3=\span(I,A_1,A_2,A_3)$.  So $A_iA_j=A_jA_i$
 for $i,j=1,2,3$.  Recall that the variety of all $(X_1,X_2,X_3)\in (\C^{4 \times 4})$ such that
 $X_iX_j=X_jX_i, i,j=1,2,3$ is irreducible \cite{GS}.  Hence a triple $(A_1,A_2,A_3)$ is a limit of
 generic commuting triples $(X_1,X_2,X_3)$.  For a generic triple, $X_1$ has $4$ distinct eigenvalues.
 So $X_2,X_3$ are polynomial in $X_1$.
 Thus, there exists $Q\in \gl(4,\C)$ such that $Q^{-1}\span (X_1,X_2,X_3)Q$ is a three dimensional subspace
 of $4\times 4$ diagonal matrices $\D\subset\C^{4\times 4}$.  Clearly $I\in\D$ and $\D$ is spanned by
 $4$ rank one diagonal matrices.  Hence $\span (X_1,X_2,X_3,I)\subseteq Q^{-1}\D Q\in\Gamma(4,\C^{4\times 4})$.
 Thus $\T_3(\cT)\subseteq\U\in \Gamma(4,\C^{4\times 4})$, and $\brank\cT\le 4$.  \qed

 In view of Lemma \ref{4invcase} we need to show Theorem \ref{landman} only in the case $\maxrank\T_p(\cT)\le 3$
 for $p=1,2,3$.   Clearly, it is enough to assume that $\cT\ne 0$.  If $\maxrank\T_p(\cT)=1$ for some $p\in\{1,2,3\}$, then $\cT_p(\cT)$ spanned by rank one matrices.  Theorem \ref{bascharrnk} implies that $\rank\cT\le 4$.
 Thus we need to consider the case
 \begin{equation}\label{uplobdrnk}
 2\le \maxrank \T_p(\cT)\le 3 \textrm{ for } p=1,2,3.
 \end{equation}

 We now consider the case $\maxrank\T_3(\cT)=2$.
 \begin{lemma}\label{countexm5.4}  Let $\cT\in\C^{4\times 4\times 4}$.  Suppose that the $(i,j)$ entry
 of each $3$-slice $T_{i,3}$ is zero if $\min(i,j)\ge 2$.  Then $\maxrank \T_3(\cT)\le 2$ and
 any three matrices $A,B,C\in\T_3(\cT)$ satisfy (\ref{commutcond1}) for $p=1,2,3$.  (In particular
 (\ref{XYZcomutcon}) holds for $\T_3(\cT)$.)  For generic choices of the four $3$-slices $T_{1,3},\ldots,T_{4,3}$ of the above form $\brank\cT=5$.  Furthermore, if
 (\ref{XYZcomutcon}) holds for $\T_1(\cT)$ and $\T_2(\cT)$, then $\brank\cT\le 4$.
 \end{lemma}
 \proof  Let $A=[a_{ij}]_{i=j=1}^4\in\C^{4\times 4}$.  Assume that $a_{ij}=0$ if $\max(i,j)\ge 2$.
 So the nonzero entries of $A$ are on the first row and column.  Clearly $\rank A\le 2$.
 Hence $\adj A=0$.  This implies that any three matrices $A,B,C\in\T_3(\cT)$ satisfy (\ref{commutcond1}) for $p=1,3$. Next observe that $\rC_2(A)$ has zero $(2,3),(2,4),(3,4)$ rows and columns.
 So $\rC_{-2}(B)$ has zero $(1,2),(1,3),(1,4)$ rows and columns.  Hence $\rC_2(A)\rC_{-2}(B)\trans\rC_2(C)=0$,
 and (\ref{commutcond1}) holds for $p=2$.

 Assume now
 \[T_{i,3}=\left[\begin{array}{cccc}a_i&b_i&c_i&d_i\\e_i&0&0&0\\f_i&0&0&0\\g_i&0&0&0\end{array}
 \right],\quad i=1,2,3,4.\]
 Consider now $T_{1,1},\ldots,T_{4,1}$.  Note that $T_{1,1}$ is a full matrix, while $T_{i,1}$ has a full first column, while the other 3 columns are equal to zero.  Suppose first that $\det T_{i,1}\ne 0$.
 (This is true of $T_{1,3},\ldots,T_{4,3}$ are generic.)
 Consider the $A_i=T_{1,1}^{-1}T_{i,1}$ for $i=1,2,3,4$.  So $A_1=I_4$ and $A_i$ has the first nonzero column
 $\a_i$, i.e. $A_i=\a_i\e_1\trans$ for $i=2,3,4$.  The commutation condition (\ref{XYZcomutcon}) with $X=A_i,Y=I_4,Z=A_j$ for $2\le i<j\le 4$ is equivalent to $(e_1\trans \a_j)\a_i\e_1\trans=(e_1\trans \a_i)
 \a_j\e_1\trans$.  Assuming that $(e_1\trans \a_j)(e_1\trans \a_i)\ne 0$ we deduce that the commutation condition holds of and only if $\a_2,\a_3,\a_4$ are colinear.  The assumption that $T_{1,3},\ldots,T_{4,3}$
 are generic matrices with nonzero entries in the first row and column yield that $\a_2,\a_3,\a_4$ are $3$
 generic vectors.  Hence $\a_2,\a_3,\a_4$ are not colinear, and the commutation condition (\ref{XYZcomutcon})
 does not hold for $\T_1(\cT)$.  Therefore $\brank\cT\ge 5$.

 To show that $\brank\cT=5$ we add to the space spanned by $I_4,\a_1\e_1\trans,\a_2\e_1\trans, \a_3\e_1\trans$ the rank one matrix $\e_1\e_1\trans$.  Let $\a_i'=\a_i-(\e_1\trans\a_i)\e_1, i=1,2,3$.
 Then the three matrices
 $A_i'=\a_i'\e_1\trans=A_i-(\e_1\trans\a_i)\e_1\e_1\trans, i=1,2,3$ commute.  Let $I_4,A_1',A_2',A_3'$ be the four $3$-slices of $\cT'\in\C^{4\times 4\times 4}$.  The proof of Lemma \ref{4invcase} yields that $\brank\cT\le 4$.  Hence $\brank\cT\le 5$ and we conclude that $\brank\cT=5$.

 We now consider non-generic $\cT\in\C^{4\times 4\times 4}$ satisfying the conditions of our lemma.
 Suppose first that $\dim\T_3(\cT)\le 2$.
 By changing a basis in the last component of $\C^4\otimes\C^4\otimes \C^4$ we may assume that
 that $T_{3,3}=T_{4,3}=0$.
 Then $\cT$ can be viewed as a tensor
 in $\C^{4\times 4\times 2}$.  (\ref{grank2mn}) yields that $\brank\cT\le 4$.
 Assume that $\dim\T_3(\cT)=3$.
 By changing a basis in the last component of $\C^4\otimes\C^4\otimes \C^4$ we may assume that
 $T_{4,3}=0$.  Assume now that the entries on the first row and column of $T_{1,3},T_{2,3}, T_{3,3}$
 are in general position. Then there exists $P=[1]\oplus P_1,Q=[1]\oplus Q_1$, where $P_1,Q_1\in\gl(3,\C)$ we may assume that $T_{i,3}'=PT_{i,3}Q=x_i\e_1\e_1\trans +\e_1\e_{i+1}\trans+\e_{i+1}\e_1 \trans$ for $i=1,2,3$.
 Let $T'_{4,3}=\e_1\e_1\trans$, and denote by $\cT'\in \C^{4\times 4\times 4}$ the tensor whose
 four $3$-slices are $T_{i,3}',i=1,2,3,4$.  We claim that $\brank\cT'\le 4$.
 Consider the following basis in $\T_3(\cT')$: $\e_1\e_1\trans$ and $\e_1\e_i\trans+\e_i\e_1\trans$ for $i=2,3,4$.
 The above arguments show that $\T_1(\cT')$ is a $4$ dimensional subspace of commuting matrices, which contain
 $I_4$.  Hence by Lemma \ref{4invcase} $\brank\cT'\le 4$.  Hence $\brank\cT\le 4$.
 Since any three matrices $T_{1,3},T_{2,3},T_{3,3}$ with zero entries in the position $(i,j)$ for $\min(i,j)\ge 2$
 can be approximated by generic matrices of this kind, we deduce that $\brank\cT\le 4$.

 We now assume that $\dim\T_3(\cT)=4$.
 Assume now that $\T_1(\cT)$ and $\T_2(\cT)$ satisfy the condition (\ref{XYZcomutcon}).
 If either $\T_1(\cT)$ or $\T_2(\cT)$ contain an invertible matrix then by Lemma \ref{4invcase}
 $\brank\cT\le 4$.  So assume that $\maxrank\T_1(\cT),\maxrank\T_2(\cT)\le 3$.
 Hence, the four first rows and columns of $T_{1,3},\ldots,T_{4,1}$ are linearly dependent.
 By changing a basis in the last component of $\C^4\otimes\C^4\otimes \C^4$ we may assume that the first
 row of $T_{4,3}$ is zero.  As $\dim \T_3(\cT)=4$ it follows that
 that $T_{4,3}$ has zero first row and nonzero first column.
 Apply elementary row operations to the last three rows of $T_{4,1}$ to assume that $T_{4,1}=\e_4\e_1\trans$.  Apply the same elementary row operations to $T_{1,3},T_{2,3},T_{3,3}$.
 Hence, we can assume without loss of generality that $T_{4,1}=\e_4\e_1\trans$.  By considering
 $T_{i,3}-t_iT_{4,3}$ for $i=1,2,3$ we may assume the $(4,1)$ entry of $T_{i,3}$ is zero for $i=1,2,3$.
 Consider again the column space spanned by the first three columns of $T_{1,3},T_{2,3},T_{3,3}$.
 It must be two dimensional, otherwise the column space $\cT$ is four dimensional and $\T_2(\cT)$
 will contain an invertible matrix.  So by changing a basis in $\span(T_{1,3},T_{2,3},T_{3,3})$
 we can assume that the first column of $T_{3,3}$ is zero.  As $\dim\T_3(\cT)=4$ the first row of $T_{3,3}$
 is nonzero.  Apply elementary column operations to the last three columns of $T_{3,3}$ we may assume that
 $T_{3,3}=\e_1\e_4\trans$.  Apply the same elementary column operations to $T_{1,3},T_{2,3},T_{4,3}$.
 We still have that $T_{4,3}=\e_4\e_1\trans$.  Apply the above arguments to deduce that without loss
 of generality we may assume that the last row and column of $T_{1,3},T_{2,3}$ are zero.
 Consider first $T_{i,2}, i=1,2,3,4$.  Observe that
 \begin{eqnarray*}
 T_{1,2}=\left[\begin{array}{cccc}*&*&0&0\\ *&*&0&0\\ *&*&0&0\\ 0&0&0&1\end{array} \right],
 T_{2,2}=\left[\begin{array}{cccc}*&*&0&0\\ 0&0&0&0\\ 0&0&0&0\\ 0&0&0&0\end{array} \right],\\
 T_{3,2}=\left[\begin{array}{cccc}*&*&0&0\\ 0&0&0&0\\ 0&0&0&0\\ 0&0&0&0\end{array} \right],
 T_{4,2}=\left[\begin{array}{cccc}0&0&1&0\\ 0&0&0&0\\ 0&0&0&0\\ 0&0&0&0\end{array} \right].
 \end{eqnarray*}
 The assumption that $T_{1,2}+T_{1,4}$ is singular yields that the second and the third row of $T_{1,2}$
 are linearly dependent.  Do elementary row operations on the second and the third row of $T_{1,2}$ to obtain that
 a zero third row.  Do the same elementary row operations on $T_{i,2}, i=2,3,4$ to deduce that we may
 assume that each $T_{i,2}$ has zero third row.  Translating back to $\T_1(\cT)$ we deduce that
 we may assume in addition to all our above assumptions on $T_{1,3},\ldots,T_{4,3}$ the third row of
 $T_{1,3},T_{2,3}$ are zero.  So all matrices in $\T_3(\cT)$ have zero third row.
 Consider now $T_{i,1}$ for $i=1,2,3,4$.  Observe that $T_{3,1}=0$.  Apply the same arguments as for $T_{i,2}$ for $i=1,2,3,4$ to deduce that we can assume that the third column of $T_{1,1}$ is zero.
 This implies that in this case we can assume that after suitable change of basis in the first two components of $\C^4\otimes\C^4\otimes\C^4$, in addition to the assumption that all $(i,j)$ entries
 of matrices in $\T_3(\cT)$ are zero if $\min(i,j)\ge 3$, the third row and column of each matrix in
 $\T_3(\cT)$ is zero.

 It is left to show that $\brank\cT\le 4$.  Observe that the last two rows and columns of
 $T_{1,3},T_{2,3}$ are zero.  Let $\cT'\in\C^{4\times 4 \times 2}$ be the tensor whose two $3$-slices are $T_{1,3},T_{2,3}$.
 So $\cT'$ can be viewed as a tensor in $\C^{2\times 2\times 2}$.  (\ref{grank2mn}) yields that $\brank\cT'\le 2$.  As $T_{3,4}=\e_1\e_4\trans, T_{4,4}=\e_4\e_1\trans$ we deduce that $\brank\cT\le 4$.
 \qed

 We remark that Lemma \ref{countexm5.4} refutes Proposition 5.4 of \cite{LM06}, which claims that
 for any tensor $\cT\in\C^{4\times 4\times 4}$ for which $\T_3(\cT)$ satisfies (\ref{XYZcomutcon})
 either $\brank\cT\le 4$ or there exists nonzero $\u\in\C^4$ such that either $T_3(\cT)\u=\0$ or
 $\u\trans\T_3(\cT)=\0\trans$.
 Indeed, if we assume as in the first part of Lemma \ref{countexm5.4} that $\cT=[t_{i,j,k}]\in\C^{4\times 4\times 4}$ is a generic tensor such that $t_{i,j,k}=0$ if $\min(i,j)\ge 2$, then there is no nonzero $\u\in\C^4$ such that either $T_3(\cT)\u=\0$ or $\u\trans\T_3(\cT)=\0\trans$, and by this lemma
 $\brank\cT=5$ .
 \begin{lemma}\label{maxranT3=2}  Let $\cT'\in\C^{4\times 4\times 4}$ and assume that $\maxrank\T_3(\cT')=2$.
 Then either $\brank\cT'\le 4$ or it is possible to change bases in the first two components of
 $\C^4\otimes\C^4\otimes\C^4$ to obtain the following two possibilities.  Either $\cT\in\C^{4\times 4\times 4}$ satisfies the conditions of Lemma \ref{countexm5.4}, or the last row and column of each matrix
 in $\T_3(\cT)$ is zero.  (In the last case, in addition every submatrix of $A\in\T_3(\cT)$ based on the first three rows and columns is singular.)
 \end{lemma}
 \proof  Change a basis in the first component of $\C^4\otimes\C^4\otimes\C^4$ such that $T_{1,3}=I_2\oplus 0$.  Since $\maxrank\T_3(\cT')=2$ it follows that the four entries $(3,3),(3,4),(4,3),(4,4)$ of each matrix in $\T_3(\cT')$ are zero.  So any matrix in $A\in\T_3(\cT')$ has the block form
 $\left[\begin{array}{cc} A_{11}&A_{12}\\A_{21}&0\end{array}\right]$, where each $A_{ij}\in\C^{2\times 2}$.
 Consider $A+tT_{1,3}$.  Assume that $\det (A_{11}+tI_2)\ne 0$.  Then $\rank (A+tT_{1,3})=2$
 if and only if $A_{21}(A_{11}+tI_2)^{-1}A_{12}=0$.
 Assume that $\rank A_{12}=2$ for some $A\in\T_3(\cT')$.  Then  for most of $A\in\T_3(\cT')$ $\rank A_{12}=2$.  Hence for most $A$'s  $A_{12}=0$.  So the last two rows of each $A\in\T_3(\cT')$ are zero.
 Hence $\cT'$ can be viewed as a tensor in $\C^{2\times 4\times 4}$ and $\brank\cT'\le 4$.
 Similarly if for $\rank A_{21}=2$ for some $A\in\T_3(\cT')$ we deduce that $\brank\cT'\le 4$.
 It is left to assume that $\rank A_{12},\rank A_{21}\le 1$ for each $A\in\T_3(\cT')$, and  $\rank A_{12}=\rank A_{21}=1$ for some matrix $A$.
 It is easy to see for $2 \times 2$ matrices that $A_{12}$ is either $\x(A_{12})\u\trans$ or $\u\x(A_{12})\trans$, and $A_{21}$ is either $\y(A_{21})\v\trans$ or $\v\y(A_{21})\trans$, where $\u,\v\in\C^2\setminus\{\0\}$ and $\x(A_{12}),\y(A_{21})\in\C^2$ are nonidentically
 zero vectors depending linearly on the entries of $A_{12}$ and $A_{21}$ respectively.
 By using elementary row or column operations on the first two rows and columns, and then on the last two rows and columns we may assume that $\u=\v=\e_1$.   We now obtain the four $3$-slices of $\cT$.
 Note that for some $T_{1,3}=B\oplus 0$, where $B\in\C^{2\times 2}$ is invertible.

 Suppose first that $A_{12}= \e_1\x(A_{12})\trans, A_{21}=\y(A_{21})\e_1\trans$.  Since $\x(A_{12})$
 and $\y(A_{21})$ are not zero identically, and $\rank A\le 2$ we deduce that the $(2,2)$ entry of each $A$ is zero.  (For example if $(3,1)$ and $(1,3)$ entries of $A$ are nonzero then consider the minor with the first three rows and columns.)
 Hence $\T_3(\cT)$ is of the form given by Lemma \ref{countexm5.4}.

 Suppose next that $A_{12}=\x(A_{12})\e_1\trans, A_{21}=\e_1\y(A_{21})\trans$.
 Then the last row and column of each matrix in $\T_3(\cT)$ is zero.
 Moreover, each  $3 \times 3$ submatrix of $\T_3(\cT)$ is singular.

 The next case is $A_{12}=\x(A_{12})\e_1\trans, A_{21}=\y(A_{21})\e_1\trans$.
 Thus the last column of each $A\in\T_3(\cT)$ is zero.
 Since $\x(A_{21})$ is not identically zero, and $\maxrank\T_2(\cT)=2$ we deduce the minor based
 on rows $1,2$ and columns $2,3$ must be identically zero.  So we have two possibilities.
 First possibility: by elementary row operations on the first two rows of matrices on $\T_3(\cT)$ we can bring $\T_3(\cT)$ to matrices of the form given by Lemma \ref{countexm5.4} with an addition condition,
 that the last column of all these matrices is zero.
 Second possibility: by elementary column operations we can achieve that also the third column
 of all matrices in $\T_3(\cT)$ are zero.  So $\cT$ can be viewed as a tensor in $\C^{4\times 2\times 4}$.
 Hence $\brank \cT'\le 4$.  Similar results hold for the last case $A_{12}=\e_1\x(A_{12})\trans, A_{21}=\e_1\y(A_{21})\trans$.  \qed

 We state the precise version of Theorem \ref{mrank=m-1} for $m=4$.
 \begin{lemma}\label{mrank=3case}  Let $\U\subset\C^{4\times 4}$ and assume that $\maxrank \U=3$.
 Then any three matrices in $\U$
 satisfy (\ref{XYZcomutcon}) if and only if one of the following mutually exclusive conditions hold.
 \begin{enumerate}
 \item\label{mrank=3a}  There exists a nonzero $\u\in\C^4$ such that either $\U\u=\0$ or $\u\trans \U=\0\trans$.
 \item\label{mrank=3char} $\dim \U=k+1\ge 2$.  There exists $P,Q\in\gl(4,\C)$ such that $P\U Q$ has a following basis $F_0,\ldots,F_k$.  The last row and column of $F_0,\ldots,F_{k-1}$ is zero, i.e. $F_i=G_i\oplus 0, F_i\in\C^{3\times 3}, i=0,\ldots,k-1$,
 $G_0=I_3$, and
 \begin{equation}\label{formG1k3}
 F_k=\left[\begin{array}{cc} G_k&\e_1\\ \e_2\trans&0\end{array}\right],\;G_k\in \C^{3\times 3}.
 \end{equation}
 Furthermore $G_k,G_1,\ldots,G_{k-1}$ have one the following possible three forms.
 \begin{enumerate}
 \item\label{itformG1}
 \begin{equation}\label{formG1}
 G_k=\left[\begin{array}{ccc}0&0&0\\0&0&0\\0&0&0\end{array}\right],\;
 G_i=\left[\begin{array}{ccc}0&a_i&b_i\\0&0&0\\0&c_i&d_i\end{array}\right], i=1,\ldots,k-1.
 \end{equation}
 \item\label{itformG2}
 \begin{equation}\label{formG2}
 G_k=\left[\begin{array}{ccc}0&0&0\\0&0&0\\h&0&0\end{array}\right],\;
 G_i=\left[\begin{array}{ccc}0&a_i&0\\0&0&0\\0&c_i&0\end{array}\right], i=1,\ldots,k-1.
 \end{equation}
 \item\label{itformG3}
 \begin{equation}\label{formG3}
 G_k=\left[\begin{array}{ccc}0&0&0\\0&0&f\\0&0&0\end{array}\right],\;
 G_i=\left[\begin{array}{ccc}0&a_i&b_i\\0&0&0\\0&0&0\end{array}\right], i=1,\ldots,k-1.
 \end{equation}
 \end{enumerate}
 \end{enumerate}
 \end{lemma}
 \proof
 By changing a basis in $\C^3$ we can assume that $\f,\g\in\C^3$ appearing in Theorem \ref{mrank=m-1} for $m=4$ are of the form $\f=\e_1,\g=\e_2$.
 Next we observe that $G_k$ can be always assumed to be of the form
 \begin{equation}\label{formG4}
 G_k=\left[\begin{array}{ccc}0&0&0\\g&0&f\\h&0&0\end{array}\right].
 \end{equation}
 Indeed, first replace $F_k$ by $F_k'=F_k-tF_0$ such that the $(3,3)$ entry of $F_k$ is zero.
 Next, use the following elementary row and column operations without changing the form of $F_0,\ldots,F_{k-1}$. Subtract a multiple of a row four from row $i$ for $i=1,2,3$.  Similarly,
 subtract a multiple of a column four from column $i$ for $i=1,2,3$.  The exact forms of $G_k,G_1,\ldots,G_{k-1}$ are obtained by choosing subspaces $X,Y$ appearing in Theorem \ref{mrank=m-1} to be of the following forms:
 $X=\span(\e_1),Y=\span(\e_2)$; $X=\span(\e_1,\e_3),Y=\span(\e_2)$; $X=\span(\e_1),Y=\span(\e_2,\e_3)$.

 Consider now the choice  $X=\span(\e_1,\e_3),Y=\span(\e_2,\e_3)$.  Then $G_i$ is of the form
 given by \eqref{formG1} where $b_i=c_i=d_i=0$ for $i=1,\ldots,k-1$.  Furthermore $G_k=0$.
 This can be considered as a special case of \ref{itformG1}.
 \qed

 \noindent
 \textbf{Proof of Theorem \ref{landman}.}  In view of Theorem \ref{bascharrnk} and Lemmas \ref{4invcase}, \ref{maxranT3=2}, \ref{countexm5.4} we need only to consider the case $\maxrank\T_p(\cT)= 3$ for $p=1,2,3$.
 Assume first that $\dim\T_p(\cT)\le 2$ for some $p\in\{1,2,3\}$.  Then by changing a basis in the $p$-th factor of $C^4\otimes\C^4 \otimes\C^4$ and permuting the factors, we obtain that $\cT$ can be viewed as a tensor in $\C^{2\times 4\times 4}$.  \eqref{grank2mn} yields that $\brank\cT\le 4$.
 Hence we assume that $\dim\T_p(\cT)\ge 3$ for $p=1,2,3$.
 Assume next that $\dim\T_p(\cT)=\dim\T_q(\cT)=3$ for some $1\le p<q\le 3$.  Then Theorem \ref{landman}
 holds.

 Assume next that $\dim\T_p(\cT)=3$ for some $p\in\{1,2,3\}$ and $\dim\T_q(\cT)=4$ for
 $q\in\{1,2,3\}\setminus\{p\}$.  By permuting the factors of $C^4\otimes\C^4 \otimes\C^4$ we may assume that
 $\dim\T_1(\cT)=\dim\T_2(\cT)=4, \dim\T_3(\cT)=3$.
 Observe first that there is no
 nonzero $\u\in\C^4$ such that either $\T_3(\cT)\u=\0$ or $\u\trans\T_3(\cT)=\0\trans$.
 Indeed, assume first that $\T_3(\cT)\u=\0$ for some nonzero $\u$.  By change of coordinates in
 the second component of $\C^4\otimes\C^4\otimes\C^4$ we can assume that $\u=\e_4$.
 So the fourth column of each matrix in $\cT_3(\cT)$ is zero.  Hence $T_{4,2}=0$ which contradicts that
 $\dim\T_2(\cT)=4$.  Similar arguments apply if $\u\trans \T_3(\cT)=\0\trans$.

 We now apply part \emph{2} of Lemma \ref{mrank=3case}.  Here $k=2$.  Assume first that $F_1,F_2$
 have the form given in \emph{2.a}.  Consider the matrix $G_1=\left[\begin{array}{ccc}0&a&b\\0&0&0\\0&c&d
 \end{array}\right]$.  Assume the generic case $d\ne 0$ and $ad-bc\ne 0$.  Then the three eigenvalues of $G_1$ are $0,0,d$.  The Jordan canonical form of $G_1$ is $J=\left[\begin{array}{ccc}0&e&0\\0&0&0\\0&0&d
 \end{array}\right]$, where $e\ne 0$.  Furthermore, $PG_1P^{-1}=J$ where  and $P,P^{-1}$ have the form
 \begin{equation}\label{exactfrmP33}
 P=\left[\begin{array}{ccc} 1&\alpha&\beta\\0&1&0\\0&\gamma&1\end{array}\right],
 \quad P^{-1}=\left[\begin{array}{ccc} 1&-\alpha+\beta\gamma&-\beta\\0&1&0\\0&-\gamma&1 \end{array}\right].
 \end{equation}
 (Note that
 $\e_1,\e_2\trans$ are right and left eigenvectors of $G_1,J$ corresponding to $0$
 eigenvalue.)  Let $Q=P\oplus [1]\in\gl(4,C)$.  Then $QF_0Q^{-1}=F_0, QF_1Q^{-1}=J\oplus[0],
 QF_2D^{-1}=F_2$.  Equivalently, we may assume that $G_1$ is equal to $J$.
 Let $F_3=\e_3\e_3\trans\in\C^{4\times 4}$ and consider the tensor $\cT'\in \C^{4\times 4\times 4}$
 whose four $3$-slices are
 \[T_{1,3}'=F_0-\e_3\e_3\trans,\; T_{2,3}'=\frac{1}{e}(F_1-d\e_3\e_3\trans),\; T_{3,3}'=F_2,\; T_{4,3}=\e_3\e_3\trans.\]
 Let $\cT''\in\C^{4\times 4\times 3}$ is obtained from $\cT'$ by deleting the last $3$-slice of $\cT'$.
 We claim that $\brank\cT''\le 3$.  Observe that the three $3$-slices of $\cT''$ have zero third row and column.  So we can view $\cT''$ as $\cS\in \C^{3\times 3\times 3}$ whose three slices are given
 as in the last part of the proof of Theorem \ref{charrnk3ten}.  Hence $\brank\cT''\le 3$
 and $\brank\cT\le 4$.

 Assume now that $F_1,F_2$ have the form given in \emph{2.c}.  So \[G_1=\left[\begin{array}{ccc}0&a&b\\0&0&0\\0&0&0
 \end{array}\right], G_2=\left[\begin{array}{ccc}0&0&0\\0&0&f\\0&0&0
 \end{array}\right].\]
 Assume the generic case that $b,f\ne 0$.  Consider
 $F_i'=(P\oplus [1])F_i(P\oplus[1])^{-1}$ for
 $i=0,1,2$.  Assume that $P$ is of the form \eqref{exactfrmP33}, where $\alpha=\beta=0, \gamma=\frac{a}{b}$.
 Then
 \[F_0'=F_0, F_1'=G_1'\oplus [0], F_2'=\left[\begin{array}{cc}G_2'&\e_1\\ \e_2\trans&0\end{array}\right],
 G_1'=\left[\begin{array}{ccc}0&0&b\\0&0&0\\0&0&0\end{array}\right],
 G_2'=\left[\begin{array}{ccc}0&0&0\\0&-\frac{fa}{b}&f\\0&-\frac{fa^2}{b^2}&\frac{fa}{b}\end{array}\right].\]
 Let $F_2''=F_2'-\frac{fa}{b}F_0$.  Then do the following elementary column and row operations on $F_0,F_1',F_2''$ to obtain $\hat F_i, i=0,1,2$.  Add to column one $\frac{fa}{b}$ times column four, add to row three $\frac{fa^2}{b^2}$ times row four and add to row two $\frac{2fa}{b}$ times row four.
 Observe that $\hat F_0,\hat F_1,\hat F_2$ are of the form $F_0,F_1,F_2$ we started with, and with the addional fact that $a=0$ in $G_1$.
 Since $b\ne 0$, by replacing $F_1$ with $\frac{1}{b}F_1$ we may assume that $b=1$.
 It is left to show that our tensor $\cT\in\C^{4\times 4\times 3}$ has border rank $4$ at most.
 Let $\cT(t)\in \C^{4\times 4\times 3}$ be the tensors with the following three $3$-slices
 $F_0(t)=I_3\oplus[t],F_1,F_2$.  It suffices to show that for any $t\ne 0$, $\brank\cT(t)\le 4$.
 Mulitply the last row of $F_0(t),F_1,F_2$ by $\frac{1}{t}$ to obtain $I_4,F_1,F_2(t)$.
 Let $\cT(t)'\in\C^{4\times 4\times 3}$ be the tensor with the above three $3$-slices.
 Observe that $F_1F_2(t)=F_2(t)F_1=0$.  Apply the arguments in the proof of Lemma \ref{4invcase}
 to deduce that $\brank\cT(t)'\le 4$.

 Assume now that $F_1,F_2$ have the form given in \emph{2.b}.
 Consider $F_0\trans=F_0,F_1\trans,F_2\trans$.  Let $P=\left[\begin{array}{ccc}0&1&0\\1&0&0\\0&0&1
 \end{array}\right]$.  Then $F_0=(P\oplus[1])F_0\trans (P\oplus[1])^{-1}$ and
 $(P\oplus[1])F_i\trans (P\oplus[1])^{-1},i=1,2$ are of the form
 \emph{2.c}.  Hence $\brank\cT\le 4$.

 It is left to discuss the case where $\dim\T_p(\cT)=4,\maxrank\T_p=3,p=1,2,3$.
 We show that this case does not exists.  As above we observe that there is no nonzero $\u$
 such that either $\T_p(\cT)\u=0$ or $\u\trans\T_p(\cT)=\0\trans$.
 We now apply Lemma \ref{mrank=3case} to $\T_3(\cT)$.  We change bases in the first two components
 of $\C^4\otimes\C^4\otimes\C^4$ to obtain one in the three possibilities discussed in Lemma \ref{mrank=3case}.
 We start with the case \emph{2.a}, where $k=3$.  Our assumption that $\dim\T_3(\cT)=4$ implies that
 the matrices $G_1$ and $G_2$ in \eqref{formG1} are linearly independent.
 Observe next that
 \begin{eqnarray*}
 T_{1,2}=\left[\begin{array}{cccc}1&0&0&0\\ 0&0&0&0\\ 0&0&0&0\\ 0&0&0&0\end{array} \right],
 T_{2,2}=\left[\begin{array}{cccc}0&a_1&a_2&0\\ 1&0&0&0\\ 0&c_1&c_2&0\\ 0&0&0&1\end{array} \right],\\
 T_{3,2}=\left[\begin{array}{cccc}0&b_1&b_2&0\\ 0&0&0&0\\ 1&d_1&d_2&0\\ 0&0&0&0\end{array} \right],
 T_{4,2}=\left[\begin{array}{cccc}0&0&0&1\\ 0&0&0&0\\ 0&0&0&0\\ 0&0&0&0\end{array} \right].
 \end{eqnarray*}
 The assumption that $\maxrank\T_2(\cT)=3$ yields that
 $\det (T_{2,2}+x_1T_{1,2}+x_3T_{3,2}+x_4T_{4,2})=0$.  Hence
 $(a_1+x_3b_1)(c_2+x_3d_2)-(a_2+x_3b_2)(c_1+x_3d_1)$ is identically zero.
 Let
 \[A_1=\left[\begin{array}{cc} a_1&a_2\\c_1&c_2\end{array}\right],\;
 A_2=\left[\begin{array}{cc} b_1&b_2\\d_1&d_2\end{array}\right].\]
 Then there exists a nonzero $\v\in\C^2$ such that either $A_1\v=A_2\v=\0$ or $\v\trans A_1=\v\trans A_2=\0\trans$.  The first possibility yields that there exists nonzero $\u\in\C^4$ such that $\T_2(\cT)\u=0$
 which contradicts our assumptions.  Hence, there exists a nonzero $\v\in\C^2$ such that  $\v\trans A_1=\v\trans A_2=\0\trans$.

 The assumption that $\maxrank\T_1(\cT)=3$ yields that
 $(a_1+x_3c_1)(b_2+x_3d_2)-(a_2+x_3c_2)(b_1+x_3d_1)$ is identically zero.
 Let
 \[B_1=\left[\begin{array}{cc} a_1&b_1\\a_2&b_2\end{array}\right],\;
 A_2=\left[\begin{array}{cc} c_1&d_1\\c_2&d_2\end{array}\right].\]
 As above we deduce that there exists nonzero $\w\in\C^2$ such that $B_1\w=B_2\w=\0$.
 In particular we deduce that the rows $(a_1,b_1),(a_2,b_2)$ and $(c_1,d_1), (c_2,d_2)$ are linearly
 dependent.  Hence, by choosing a new basis in $\span(G_1,G_2)$ we can assume that $c_1=d_1=0$
 and $a_2=b_2=0$.  So $A_1,A_2$ are both diagonal, singular, and have a common left zero eigenvector.
 So either $a_1=b_1=0$ or $c_2=d_2=0$.  That is either $G_1=0$ or $G_2=0$ which contradicts our assumption that $\dim \T_3(\cT)=4$.

 Assume the condition \emph{2.b}.  Consider the matrices $T_{1,2},T_{2,2},T_{3,2},T_{4,2}$.
 Note that $ T_{2,2}=\left[\begin{array}{cccc}0&a_1&a_2&0\\ 1&0&0&0\\ 0&c_1&c_2&0\\ 0&0&0&1\end{array} \right]$.
 The assumption that $\maxrank\T_2(\cT)=3$ yields that $a_1c_2-a_2c_1=0$.  So the vectors $(a_1,c_1)\trans, (a_2,c_2)\trans$ are linearly dependent.  Hence $F_1,F_2$ are linearly dependent.  This contradicts
 the assumption that $\dim\T_3(\cT)=4$.

 Assume the condition \emph{2.c}.  Consider the matrices $T_{1,1},T_{2,1},T_{3,1},T_{4,1}$.
 The assumption that $\det T_{1,1}=0$ yields that the vectors $(a_1,b_1), (a_2,b_2)$ are linearly
 dependent.  Hence $F_1,F_2$ are linearly dependent.  This contradicts
 the assumption that $\dim\T_3(\cT)=4$.
 \qed

 \noindent
 \textbf{Proof of Theorem \ref{charbr4}.}  We first show that if $\brank\cT\le 4$ then the conditions \ref{charbr4a}-\ref{charbr4b}
 hold.  Theorem \ref{commutcond} yields the condition \ref{charbr4a}.  Let $P_1,P_2,P_3\in \gl(4,\C)$.
 Then $\brank \cT(P_1,P_2,P_3)\le 4$.  Let $\cT'(P_1,P_2,P_3)=[t_{i,j,k}(P_1,P_2,P_3)]_{i=j=k}^{3,3,4}\in \C^{3\time 3\times 4}$.  Clearly $\brank\cT'(P_1,P_2,P_3)\le 4$.  Theorem \ref{334charbr4} yields the conditions (\ref{rankcondTP123a})-(\ref{rankcondTP123b}) for $p=3$ and (\ref{symmet1}) for $R_3(P_1,P_2,P_3),
 L_3(P_1,P_2,P_3)$.  Similar arguments imply the conditions \ref{charbr4b} for $p=1,2$.  Since $\gl(4,\C)$ is dense in $C^{4\times 4}$ we deduce the condition \ref{charbr4b} for any
 $P_1,P_2,P_3\in \C^{4\times 4}$.

 Assume to that $\cT$ satisfies the conditions \ref{charbr4a}-\ref{charbr4b}.  Suppose to the contrary that
 $\brank\cT>4$.  Theorem \ref{landman} yields that there exists nonzero $\u\in\C^4$ such that either
 $\T_3(\cT)\u=\0$ or $\u\trans\T_3(\cT)=0$.  By changing the first two factors if needed we can assume
 that $\u\trans\T_3(\cT)=\0$.  By changing a basis in the first factor we may assume that $\u=\e_4=(0,0,0,1)\trans$.  Consider now $\T_1(\cT)$.  Note that the fourth $1$-slice $T_{4,1}$ is zero
 matrix.  Since $\brank\cT>4$ Theorem \ref{landman} yields that there exists nonzero $\v\in\C^4$ such that
 either $\T_1(\cT)\v=\0$ or $\v\trans\T_1(\cT)=\0$.  By permuting the two last factors $\C^4$, if needed, we may assume that $\v\T_1(\cT)=\0$.  By changing a basis in the second factor we can assume that $\v=\e_4$.
 This finally means that after permuting the factors of $\C^4\times \C^4\times \C^4$, and changing bases
 in the first two factors, we obtain a new tensor $\cT'$ whose four $3$-slices are matrices with the zero last row and column.  Permuting back to the original factors and changing bases correspondingly using $P_1,P_2,P_3\in \gl(4,\C)$, (one of these matrices is an identity matrix), we deduce that
 for some $p\in\{1,2,3\}$, the four $p$-slices $S_{1,p}(P_1,P_2,P_3),\ldots,S_{4,p}(P_1,P_2,P_3)$
 have zero last row and column.  Combine conditions \ref{charbr4b} with Theorem \ref{334charbr4}
 to deduce that
 the tensor $\cT'\in\C^{3\times 3\times 4}$, whose four $3$-slices are $T_{1,p}(P_1,P_2,P_3),\ldots,T_{4,p}(P_1,P_2,P_3)$, has border rank $4$ at most.
 Hence $4\ge \brank\cT(P_1,P_2,P_3)=\brank \cT$ which contradicts our assumption.  \qed

 We conclude this section by showing that the conditions in of Theorem \ref{charbr4} can be stated as a finite number of polynomial equations in degrees $5,9$ and $16$ in the entries of $\cT$.  First we consider the conditions \ref{charbr4a}.
 In the notation of \ref{charbr4a} $\T_p(\cT)$ is spanned by $S_{1,p}(I,I,I),\ldots,S_{4,p}(I,I,I)$.
 Fix $p\in\{1,2,3\}$.  Let
 \[X=\sum_{i=1}^4 x_iS_{i,p}(I,I,I),\quad Y=\sum_{j=1}^4 y_j S_{j,p}(I,I,I), \quad Z=\sum_{k=1}^4 z_k S_{k,p}(I,I,I).\]
 Then $\adj Y$ is $4\times 4$ matrix whose entries are homogeneous polynomials
 of degree $3$ in $ \y=(y_1,\ldots,y_4)\trans$.
 Note that the  coefficients of the monomials of these polynomials are polynomials of degree
 $3$ in the entries of $\cT$.  Substitute the expressions of $X,\adj Y,Z$ into the conditions
 (\ref{XYZcomutcon}) to deduce that a finite number of polynomials in $\x,\y,\z$ of degree $5$
 must vanish identically.  That is, the corresponding coefficient of each monomial
 must be zero.  This procedure gives rise to a finite number of polynomial equations of degree $5$ that the entries of $\cT$ satisfy.  Clearly, if these conditions hold then the condition \ref{charbr4a} hold.

 We now discuss the conditions \ref{charbr4b}.  Write $P_q=[x_{ij,q}]_{i=j=1}^4$ for $q=1,2,3$.
 We view each $x_{ij,q}$ as a variable.  So the entries of $P_1,P_2,P_3$ give rise to $48$ variables.
 The entries $S_{1,p}(P_1,P_2,P_3),\ldots,S_{4,p}(P_1,P_2,P_3), p=1,2,3,$ are multilinear polynomials
 of degree $3$ in $48$ variables.  Hence the entries of $T_{1,p}(P_1,P_2,P_3),\ldots,T_{4,p}(P_1,P_2,P_3)$,
 are multilinear polynomials of degree $3$, whose coefficients are homogeneous polynomials of degree $1$
 in the entries of $\cT$.  The conditions (\ref{rankcondTP123a}-\ref{rankcondTP123b}) yield
 that each $9\times 9$ minor of the two matrices in (\ref{rankcondTP123a}-\ref{rankcondTP123b})
 must be identically zero.  Such a minor is a homogeneous polynomial of degree$27=3\times 9$ in the
 entries of $P_1,P_2,P_3$.  The coefficient of each monomial is a polynomial of degree $9$ in the entries
 of $\cT$.  Hence the coefficient of each monomial appearing in the expansion of each minor must equal to zero.  (Recall that for each $p\in\{1,2,3\}$ we have $440$ such $9\times 9$ minors.)  These conditions give rise to a finite number polynomial conditions of degree $9$ on the entries of $\cT$.
 Repeat the same procedure for the conditions on $R_p(P_1,P_2,P_3), L_p(P_1,P_2,P_3)$ given by (\ref{sufconsym1}-\ref{sufconsym2}) to deduce a finite number of polynomial conditions of degree $16$
 on the entries of $\cT$.  Clearly, these polynomial equations of degree $9$ and $16$ imply that
 the condition \ref{charbr4b} hold for each $P_1,P_2,P_3\in\C^{4\times 4\times 4}$.\\

 \section{Tensors in $\C^{m\times n\times l}$ of rank $l$}
 Let $m,n,l\ge 2$ and and assume that $\cT=[t_{i,j,k}]\in \C^{m\times n\times l}$.
 In this section we study mainly the conditions
 when $\rank \cT= \dim\T_3(\cT)$.
 We point out briefly how to state some of these conditions for tensors of border rank $l$ at most.
 We consider the generic case $\dim\T_3(\cT)=l$.
 Equivalently we study the conditions on a subspace $\T\subset \C^{m\times n}$,
 of dimension $l$, to be spanned by $l$ linearly independent rank one matrices.
 First choose a basis $T_1,\ldots,T_l\in\C^{m\times n}$ of $\T$.
 Let $\z=(z_1,\ldots,z_l)\trans\in\C^l$ and denote $T(\z)=\sum_{k=1}^l z_k T_k$.
 Recall that $T(\z)$ has rank at most rank 1 if all $2\times 2$ minors of $T(\z)$ are zero.
 Let
 $$\alpha=(i_1,i_2), 1\le i_1<i_2\le m, \quad \beta=(j_1,j_2), 1\le j_1<j_2\le n.$$
 As in \S2 denote by $2_2^{\an{m}},2_2^{\an{n}}$ the set of all allowable $\alpha,\beta$ respectively.
 Then $T(\z)[\alpha,\beta]$ is the $2\times 2$ minor of $T(\z)$ based on the rows $i_1,i_2$
 and the columns $j_1,j_2$.  Clearly $T(\z)[\alpha,\beta]$ is a quadratic form in $\z$.
 For $A=[a_{ij}],B=[b_{ij}]\in \C^{m\times n}$ denote
 $$b(A,B)[\alpha,\beta]=\det\left[\begin{array}{cc}a_{i_1j_1}&b_{i_1j_2}\\
 a_{i_2j_1}&b_{i_2j_2}\end{array}\right].$$
 Then the condition that $T(\z)$ has rank one at most is given by the system
 of quadratic equations
 \begin{equation}\label{secondmin0}
 T(\z)[\alpha,\beta]=\sum_{p,q=1}^l b(A_p,A_q)[\alpha,\beta] z_pz_q=0, \; \textrm{ for all }
 \alpha\in 2_2^{\an{m}}, \beta\in 2_2^{\an{n}}.
 \end{equation}
 Note the number of these equations is ${m \choose 2}{n \choose 2}$.
 With each of the above quadratic form we associate a symmetric matrix
 $S(\alpha,\beta)(T_1,\ldots,T_l)\in \rS(l,\C)$,
 Denote by $\bS(T_1,\ldots,T_l)\subset \rS(l\C)$ the subspace spanned by
 $S(\alpha,\beta)(T_1,\ldots,T_l),\alpha\in 2_2^{\an{m}}, 2_2^{\an{n}}$.
 Suppose that we change a basis of $\T$ from $T_1,\ldots,T_l$ to $T_1',\ldots,
 T_l'$.  This is equivalent to the change of variables $\z=R\y$.  Hence
 \begin{equation}\label{changbas}
 S(\alpha,\beta)(T_1',\ldots,T_l')=R\trans S(\alpha,\beta)(T_1,\ldots,T_l)R,\textrm{ for some }
 R\in\gl(l,\C),
 \end{equation}
 and  $\alpha\in 2_2^{\an{m}},\beta\in 2_2^{\an{n}}$.  Thus $\bS(T_1',\ldots,T_l')=R\trans \bS(T_1,\ldots,T_l)R$.
 For simplicity of notation we let $\bS(\T):=\bS(T_1,\ldots,T_l)$.  So the subspace $\bS(\T)$ is defined
 up to congruence.
 \begin{prop}\label{PTQprop}  Let $T_1,\ldots,T_l$ be a basis in $\T\subset\C^{m\times n}$.
 Then for any $P\in \gl(m,\C),Q\in\gl(n,\C)$
 \begin{equation}\label{PTQprop1}
 S(\alpha,\beta)(PT_1Q,\ldots,PT_lQ)=\sum_{\gamma\in 2_2^{\an{m}},\delta\in 2_2^{\an{n}}}
 P[\alpha,\gamma]Q[\delta,\beta] S(\gamma,\delta)(T_1,\ldots,T_l).
 \end{equation}
 Hence $\bS(P\T Q)=\bS(\T)$.
 \end{prop}
 \proof
 Recall that for any $A\in \C^{m\times n}$ the $2 \times 2$ minors of $PAQ$ are given by the Cauchy-Binet
 formula
 $$(PAQ)[\alpha,\beta]=\sum_{\gamma\in 2_2^{\an{m}},\delta\in 2_2^{\an{n}}} P[\alpha,\gamma]A[\gamma,\delta]
 Q[\delta,\beta].$$
 Clearly $PT(\z)Q=\sum_{k=1}^l z_k PT_kQ$.  Apply the Cauchy-Binet formula to deduce (\ref{PTQprop1}).
 Hence $\bS(P\T Q)\subseteq \bS(\T)$.  Since $\T=P^{-1}(P\T Q)Q^{-1}$ we obtain that
 $\bS(P\T Q)\supseteq \bS(\T)$.
 \qed

 Let $\bS(T_1,\ldots,T_l)^{\perp}\subset \rS(l,\C)$ be the orthogonal complement
 of with respect the symmetric product on $\C^{l\times l}$: $\an{A,B}:=\tr AB\trans$.
 \begin{lemma}\label{rank1perlem}  Let $\T\subset \C^{m\times n}$ be an $l$-dimensional subspace.
 Then $\T$ contains $r$-linearly independent rank one matrices if and only the subspace $\bS(\T)^{\perp}$
 contains $r$ linearly independent rank one symmetric matrices which are simultaneously diagonable by
 congruency.  That is, if $T_1,\ldots,T_l$ is a basis of $\T$, then there exist $R\in\gl(l,\C)$
 such that

 \noindent
 $R\trans \diag(\delta_{1k},\ldots,\delta_{lk})R\in\bS(\T)^\perp$ for $k=1,\ldots,r$.
 \end{lemma}
 \proof  Choose a basis in $\T$: $T_1,\ldots,T_l$ such that $T_1,\ldots,T_r$ are $r$-rank one linearly
 independent matrices.  Hence the quadratic form $T(\z)[\alpha,\beta]$ does not contain terms
 $z_1^2,\dots,z_r^2$.  Thus the diagonal entries $(1,1),\ldots,(r,r)$ are zero for each
 $S(\alpha,\beta)(T_1,\ldots,T_l)$.  Therefore $D_k:=\diag(\delta_{1k},\ldots,\delta_{lk})\in
 \bS(T_1,\ldots,T_l)^{\perp}$ for $k=1,\ldots,r$.

 Assume now that $\bS(T_1,\ldots,T_k)^{\perp}$ contains $R^{-1}D_k(R^{-1})^{\trans}$ for $k=1,\ldots,r$.
 This is equivalent to the fact that $D_k\in \bS(T_1',\ldots,T_l')^{\perp}$ for a corresponding basis
 $T_1',\ldots,T_l'$ of $\T$.  Hence $T_1',\ldots,T_r'$ are $r$ linearly independent rank one matrices
 in $\T$.  \qed

 \begin{corol}\label{dimest}  Let $\T\subset \C^{m\times n}$ be an $l$-dimensional subspace.
 Assume that $\T$ contains $r$-linearly independent rank one matrices.
 Then dimension $\bS(\T)$ is at most ${l+1 \choose 2} -r$.
 \end{corol}
 The data induced by  $S(\alpha,\beta)(T_1,\ldots,T_l)$ for $\alpha\in 2_2^{\an{m}},\beta\in 2_2^{\an{n}}$
 can be arranged in the following ${m \choose 2}{n \choose 2} \times {l+1 \choose 2}$ matrix
 $C(\T)=C(T_1,\ldots,T_l)=[c_{(\alpha,\beta)(p,q)}]$.  Here
 \begin{equation}\label{defcalphabetapq}
 c_{(\alpha,\beta)(p,q)}=b(T_p,T_q)[\alpha,\beta]+b(T_q,T_p)[\alpha,\beta],\;1\le p\le q\le l,
 \alpha\in 2_2^{\an{m}},\beta\in 2_2^{\an{n}}.
 \end{equation}
 Corollary \ref{dimest} equivalent to the statement that $\rank C(\T)\le {l+1\choose 2}-r$.
 That is, all minors of order  ${l+1\choose 2}-r+1$ of $C(\T)$ are zero.
 We now give two examples of generic subspaces $\T\subset \C^{m\times n}$ of dimension $l$ spanned by rank one matrices which satisfy $\dim\rS(\T)^\perp = l$.
 In these case we obtain necessary conditions for $\cT\in\C^{m\times n\times l}$ to have a border rank
 $l$ at most.
 \begin{lemma}\label{genmngelcase}
 Assume that one of the following conditions hold.
 \begin{enumerate}
 \item\label{genmngelcase1}
 $2\le l\le m,n$ and $\T\subset \C^{m\times n}$ is an $l$-dimensional
 subspace spanned by $l$ rank one matrices $\u_1\v_1\trans,\ldots,\u_l\v_l\trans$, where $\u_1,\ldots,\u_l
 \in\C^m$ and $\v_1,\ldots,\v_l\in\C^n$ are linearly independent.
 \item\label{genmngelcase2} $m=n=l-1\ge 3$ and $\T\subset \C^{(l-1)\times (l-1)}$
 is an $l$-dimensional subspace spanned by $l$ rank one matrices $\u_1\v_1\trans,\ldots,\u_l\v_l\trans$, where any $l-1$ vectors out of $\u_1,\ldots,\u_l
 \in\C^{l-1}$ and $\v_1,\ldots,\v_l\in\C^{l-1}$ are linearly independent.
 \end{enumerate}
 Then $\dim \bS(\T)={l\choose 2},\dim\bS(\T)^\perp=l$, and $\bS(\T)^\perp$
 spanned by $l$ ranks one symmetric matrices which are simultaneously diagonable by congruency.
 Hence for any $\T\in\Gamma(l,\C^{m\times n})$ the matrix $C(\T)$ has rank at most
 ${l\choose 2}$.  In particular if $\cT=[t_{i,j,k}]\in\C^{m\times n\times l}$ has border rank not more than $l$,
 then the $l$ $3$-slices $T_k:=[t_{i,j,k}]_{i=j=1}^{m,n}\in\C^{m\times n}, k=1,\ldots,l$ satisfy the identities
 given by the vanishing of all ${l\choose 2}+1$ minors of $C(T_1,\ldots,T_l)$.
 \end{lemma}
 \proof  Let $\e_i:=(\delta_{i1},\ldots,\delta_{im})\trans, \f_j:=
 (\delta_{j1},\ldots,\delta_{jn})\trans$ for $i=1,\ldots,m$ and $j=1,\ldots,n$.
 Consider first that the case \ref{genmngelcase1}.
 Since $\u_1,\ldots,\u_l$ and $\v_1,\ldots,\v_l$ are linearly independent, it follows that there exist
 $P\in \gl(m,\C), Q\in\gl(n,\C)$ such that $P\u_i=\e_i, Q\v_i=\f_i, i=1,\ldots,l$.
 Let $T_i':=P(\u_i\v_i\trans)Q\trans=\e_i\f_i\trans$ for $i=1,\ldots,l$.
 Thus $T'(\z)=\diag(z_1,\ldots,z_l,0,\ldots,)$.  The only nonzero $2\times 2$ minors of $T'(\z)$
 are $T'(\z)[\alpha,\alpha]$ where $\alpha=(p,q)$ and $1\le p<q\le l$.
 So $\bS(PT_1 Q\trans,\ldots,PT_lQ\trans)$ consists
 of all symmetric matrices $A\in\rS(l,\C)$ with
 zero diagonal.  Hence $\dim\bS(\T)={l \choose 2},\dim\bS(\T)^\perp=l$ and $\bS(\T)^\perp$
 spanned by $l$ ranks one symmetric matrices which are simultaneously diagonable by congruency.

 Consider now the case \ref{genmngelcase2}.  The arguments of the proof of Lemma \ref{symmet2}
 yield that we may assume that $\u_i=\v_i=\e_i$ for $i=1,\ldots,l-1$ and $\u_l=\v_l=
 \w=(w_1,\ldots,w_{l-1})\trans$, where $w_i\ne 0$ for $i=1,\ldots,l-1$.
 Let $T_i=\u_i\v_i\trans$ for $i=1,\ldots,l$.  A straightforward calculation shows that $\det T(\z)[(1,3),(1,2)]=
 w_2w_3z_1z_l$. Similarly, we have all the quadratic forms of the form $w_jw_kz_iz_l$ for $i=2,\ldots,l-1$, were $i,j,k$ are three distinct elements in $\an{l-1}$.
 By considering $\det T(\z)[\alpha,\alpha], \alpha\in 2_2^{\an{l-1}}$ we deduce that the space of quadratic
 polynomials spanned by $\det T(\z)[\alpha,\beta], \alpha,\beta\in 2_2^{\an{l-1}}$ contains all the monomials
 $z_iz_j$ for $i\ne j$.  Hence $\bS(\T)$ is all $A\in\rS(l,\C)$ with zero diagonal entries.  I.e.
 $\dim\bS(\T)={l \choose 2},\dim\bS(\T)^\perp=l$ and $\bS(\T)^\perp$
 spanned by $l$ ranks one symmetric matrices which are simultaneously diagonable by congruency.
 Other claims of the lemma follow straightforward from the continuity argument. \qed
 \begin{theo}\label{necsufrnkl} Assume that $\cT\in\C^{m\times n\times l}$,
 $\dim\T_3(\cT)=l$ and $\dim\bS(\T_3(\cT))^\perp=l$.  Then
 \begin{enumerate}
 \item\label{necsufrnkl1} $\rank\cT=l$ if and only if the following conditions hold
 \begin{enumerate}
 \item\label{necsufrnkl1a}
 $\bS(\T_3(\cT))^\perp$ contains an invertible matrix.
 \item\label{necsufrnkl1b}
 The condition (\ref{XYZcomutcon}) holds for any three matrices $A,B,C\in\bS(\T_3(\cT))^\perp$.
 \item\label{necsufrnkl1c}
 $A\;\adj B$ have $l$ distinct eigenvalues for some $A,B\in\bS(\T_3(\cT))^\perp$.
 \end{enumerate}
 \item\label{necsufrnkl2}   Assume that either $2\le l\le m,n$ or $m=n=l-1\ge 3$.
 If $\brank\cT=l$ and $\dim\bS(\T_3(\cT))^\perp=l$ then  $\bS(\T_3(\cT))^\perp\in \Gamma(l,\rS(l,\C))$.
 In particular, the conditions (\ref{commutcond1}) holds for $p=1,\ldots,l-1$ and
 any three matrices $A,B,C\in\bS(\T_3(\cT))^\perp$.
 \end{enumerate}
 \end{theo}
 \proof Let $\U:=\bS(\T_3(\cT))^\perp\subset \rS(l,\C)$.
 Suppose that $\rank\cT=l$. Lemma \ref{rank1perlem} yields that $\U$
 contains $l$ linearly independent rank one symmetric matrices which are simultaneously diagonable by a
 congruency.  Since $\dim\U=l$ we deduce that $\U$
 has a basis of the form $R\trans \diag(\delta_{1k},\ldots,\delta_{lk})R\in\rS(\T)^\perp$ for $k=1,\ldots,l$.
 Hence the conditions \ref{necsufrnkl1a}-\ref{necsufrnkl1c} hold.

 Suppose that the condition \ref{necsufrnkl1a} holds.  By considering $\cT(I_m,I_n,R)$
 for some $R\in\gl(l,\C)$ we may assume that $I\in \U$.
 The condition \ref{necsufrnkl1b} for $B=I$ yields that $\U$ is a subspace
 of commuting matrices.  Hence, there exists a unitary matrix $V$ such that $V\U V^*$ is an upper
 triangular matrix, e.g. \cite[\S24.2, Fact 3]{Hog}.
 The assumption \ref{necsufrnkl1c} yields that most of the matrices
 of the the form $A\adj B$ have simple eigenvalues. Choose $B\in \gl(l,\C)\cap\U$.
 So $\adj B = \frac{1}{\det B} B^{-1}$.  Thus, most of the matrices of the form $AB^{-1}$ have simple
 eigenvalues. Hence most of the matrices in $\U$ have simple eigenvalues.  Choose $A\in\U$
 such that $A$ has simple eigenvalues.  Thus any $C\in\U$ is a polynomial in $A$.  Since $\dim\U=l$
 it follows that $\U$ has a basis of $l$-rank one commuting matrices which are simultaneously
 diagonable by an orthogonal matrix.  Lemma \ref{rank1perlem} yields that $\rank \cT=l$.

 Assume that $\brank\cT=l$.  Hence $\cT$ is a limit of rank $l$
 tensors $\cT_q, q\in\N$ satisfying the conditions of Lemma \ref{genmngelcase}.
 Since $\dim \bS(\T_3(\cT))^\perp =l$ we deduce that $\bS(\T_3(\cT))^\perp$ is the limit
 of $\bS(\T_3(\cT_q))^\perp, q\in\N$.  (Without this assumption we can only deduce that
 any convergent subsequence of subspaces in the sequence $\bS(\T_3(\cT_q))^\perp, q\in\N$
 converges to a subspace of $\bS(\T_3(\cT))^\perp$.)
 Apply for each $\cT_q$ part \ref{necsufrnkl1} to deduce part \ref{necsufrnkl2}.  \qed

  Note that the a simultaneous matrix diagonalization by congruence arises naturally in finding the rank
 decomposition of tensors \cite{DeL06}.
 As in our characterization of $V_4(3,3,4)$ we can restate the conditions (\ref{necsufrnkl1b})
 of Theorem \ref{necsufrnkl} in terms of some polynomial equations.  These equations will also hold
 for any $\cT$ satisfying $\brank\cT\le l$.

 \bibliographystyle{plain}

 \emph{Acknowledgement}:  I thank J. M. Landsberg and the referees
 for useful remarks.

\end{document}